\documentstyle[amsfonts,amssymb,graphicx]{article}

\newtheorem{theorem}{Theorem}[section]
\newtheorem{lemma}[theorem]{Lemma}

\newtheorem{corollary}[theorem]{Corollary}

\newtheorem{proposition}[theorem]{Proposition}

\title{\bf INTEGRAL TRANSFORM APPROACH TO SOLVING KLEIN-GORDON EQUATION WITH VARIABLE COEFFICIENTS }

\author{{\bf Karen Yagdjian} }

\begin{document}

\date{}
\maketitle
\thispagestyle{empty}
\vspace{-0.3cm}

\begin{center}
{\it Department of Mathematics\\
University of Texas-Pan American\\
1201 W.~University Drive\\
Edinburg, TX 78539
USA }
\end{center}

\addtocounter{section}{-1}
\renewcommand{\theequation}{\thesection.\arabic{equation}}
\setcounter{equation}{0}
\pagenumbering{arabic}
\setcounter{page}{1}
\thispagestyle{empty}

\hspace{2cm}\begin{abstract}
\begin{small}
In this paper we describe the integral transform that allows to write solutions of  one partial differential equation via solution of another one.
This transform was suggested by the author in the case when the last equation is a wave equation, and then it was used
to investigate several well-known equations such as Tricomi-type equation, the Klein-Gordon equation in
the de~Sitter and Einstein-de~Sitter spacetimes.  A generalization given in this paper allows us to consider also
the  Klein-Gordon equations with coefficients  depending on the spatial variables.\\


\end{small}
\end{abstract}

\section{Introduction and Statement of Results}
\label{S1}

\setcounter{equation}{0}
\renewcommand{\theequation}{\thesection.\arabic{equation}}

In this paper we give some generalization of the  approach suggested in   \cite{YagTricomi}, which is aimed  to reduce  equations with variable coefficients to more simple ones.
This transform was used in a series of papers \cite{Galstian-Kinoshita-Yagdjian,Galstian-Kinoshita},
\cite{YagTricomi}-\cite{JMP2013}
to investigate in a unified way several equations such as the linear
and semilinear Tricomi  
equations, Gellerstedt
equation, the wave equation in Einstein-de~Sitter  spacetime, the
wave and the Klein-Gordon equations in the de~Sitter and
anti-de~Sitter spacetimes.  The listed  equations play an important
role in the gas dynamics, elementary particle physics, quantum field
theory in curved spaces, and cosmology.
\smallskip

Consider for the smooth function $f=f(x,t)$
the solution $w=w_{A,f}(x,t;b)$ to the  problem
\begin{equation}
\label{0.1JDE}
v_{tt}-  A(x,\partial_x) v =0, \,\, \,
v(x,0;b)= f (x,b), \,\, v_t(x,0)=0 ,\,\, t \in [ 0,T_1]\subseteq {\mathbb R} , \,\, x \in \Omega \subseteq {\mathbb R}^n,
\end{equation}
with the parameter $b\in I =[t_0,T] \subseteq {\mathbb R}$, $t_0< T \leq \infty$,  and with $ 0< T_1 \leq \infty $. Here $\Omega  $ is a domain in ${\mathbb R}^n $, while  $A(x,\partial_x) $ is the partial differential operator $A(x,\partial_x)=\sum_{|\alpha| \leq m } a_\alpha (x)D_x^\alpha $.
We are going to present   the integral operator
\begin{equation}
\label{KJDE}
{\mathcal K} [w](x,t)
 =
2   \int_{  t_0}^{t} db
  \int_{ 0}^{ |\phi (t)- \phi (b)|}   K (t;r,b;M)  w(x,r;b )  dr
, \quad x \in \Omega , \,\, t \in I,
\end{equation}
which maps the function $w=w(x,r;b ) $ into solution of the equation
\begin{equation}
\label{0.3JDE}
u_{tt}-  a^2(t) A(x,\partial_x) u -M^2 u =f, \qquad x \in \Omega \,, 
 \,\,   t \in I.
\end{equation}
In fact, the function $u=u(x,t) $ takes initial values as follows
\[
u(x,t_0)= 0, \,\, u_t(x,t_0)=0 ,\quad  x \in \Omega \,.
\]
Here $\phi =\phi (t) $ is a distance function produced by $a=a(t) $, that is $\ \phi (t) = \int_{t_0}^t a(\tau )\,d\tau $, while $M \in{\mathbb C}$ is a constant.
Moreover, we also give the corresponding operators, which generate solutions of the source-free equation and takes non-vanishing initial values. These operators are constructed in
\cite{Yag_Galst_CMP,yagdjian_DCDS} in the case of $A(x,\partial_x)=\Delta  $, where $\Delta   $ is the Laplace operator on ${\mathbb R}^n $, and, consequently, the equation (\ref{0.1JDE}) is the wave equation.
In the present paper we restrict ourselves to the smooth functions, but it is evident that  similar formulas, with the corresponding interpretations, are applicable to the distributions  
as well.  (For details see, e.g., \cite{Yag_Galst_CMP}.)

In order to motivate our approach, we consider
the solution $v=v(x,t;b)$ to the Cauchy problem
\begin{equation}
\label{1.7new}
v_{tt}-  \Delta  v =0, \,\,   (t,x) \in {\mathbb R}^{1+n},  \,\,
v(x,0;b)= \varphi  (x,b), \,\, v_t(x,0)=0 ,\,\, x \in {\mathbb R}^n,
\end{equation}
with the parameter $b \in I \subseteq {\mathbb R}$. We denote that solution by $v_\varphi =v_\varphi (x,t;b)$;
if $\varphi $ is independent of the second time variable $b$, then
  we  write simply $v_\varphi (x,t)$.
There are well-known explicit representation formulas for the solution of the   problem (\ref{1.7new}). (See, e.g.,  
\cite{Strauss}.)

The starting point of the approach suggested in \cite{YagTricomi} is the Duhamel's principle (see, e.g.,  
\cite{Strauss}), which has been revised in order to prepare the ground for generalization.
Our {\it first observation} is that the function
\begin{equation}
\label{main}
u(x,t)= \int_{t_0}^t \,db \int_{ 0}^{  t-b  } w_f(x,r;b )\,dr\,,
\end{equation}
is the solution of the Cauchy problem
$
 u_{tt}-\Delta u =f(x,t)$ in ${\mathbb R}^{n+1}$, and
$u(x,t_0)=0$, $u_t(x,t_0)=0$   in  $\,\, {\mathbb R}^{n}\,,
$
if the function $w_f=w_f(x;t;b ) $ is a solution of the problem  (\ref{1.7new}), where $\varphi =f $.
The  {\it second observation} is that in (\ref{main}) the upper limit $t-b $  of the inner integral is generated by the  propagation phenomena  with the speed which equals to one.
In fact, that is a distance function.
Our {\it third observation}  is that the solution operator ${\mathcal G}\,:\,f \longmapsto u $ can be regarded as a composition of two operators. The first one
\[
{\mathcal W}{\mathcal E}: \,\, f  \longmapsto  w
\]
is a Fourier Integral Operator, which is a solution operator of the Cauchy problem
for wave equation.  The second operator
\[
{\mathcal K}:\,\,w\longmapsto   u
\]
is the integral operator given by (\ref{main}). We regard the variable $b$ in  (\ref{main}) as a ``subsidiary time''. Thus, ${\mathcal G}= {\mathcal K}\circ {\mathcal W}{\mathcal E}$ and we arrive at the  diagram of Figure 1.
\begin{center}
\begin{figure}[hi]
\hspace{1.8cm}\includegraphics[width=0.25\textwidth]{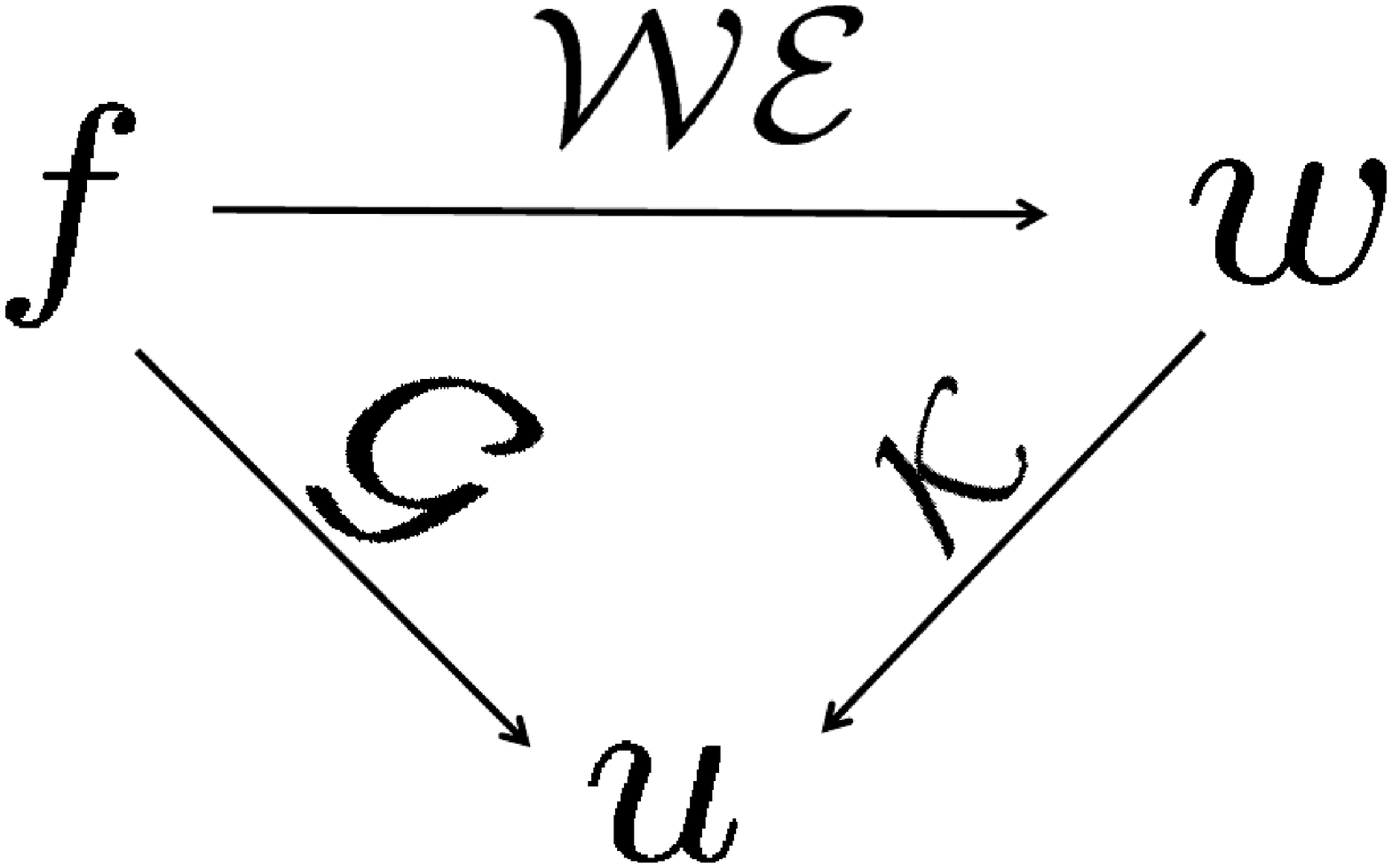}\hspace{2.9cm}
\includegraphics[width=0.25\textwidth]{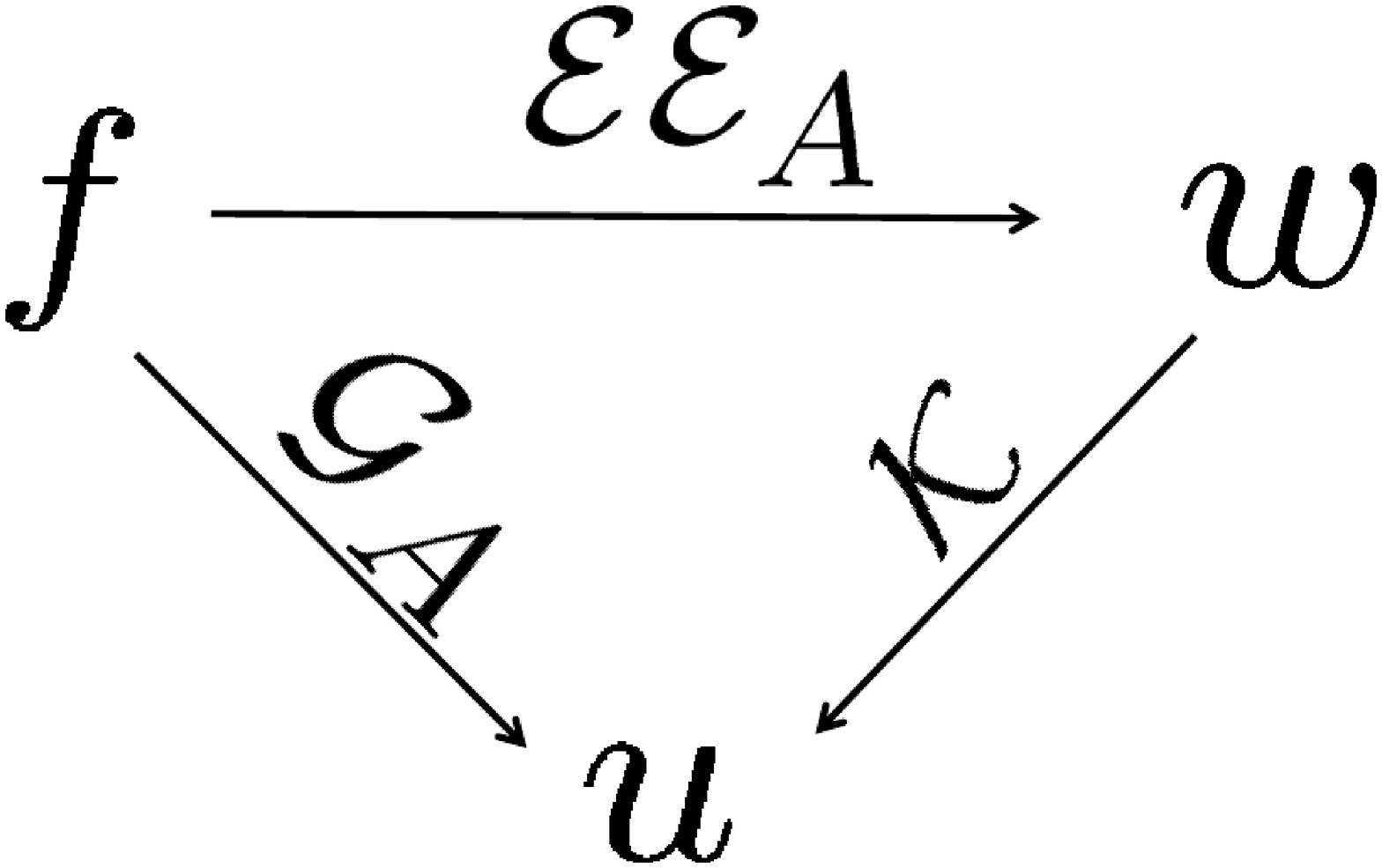}
\caption{\small Case of wave equation $A(x,\partial_x)=\Delta  $ \hspace{0.4cm} { Figure 2:}\, Case of general  $A(x,\partial_x)$}
\end{figure}
\end{center}
\vspace{-0.5cm}
\noindent
Based on the first diagram, we have generated in \cite{yagdjian_Rend_Trieste} a class of operators for which we have  obtained   explicit representation formulas for the solutions, and, in particular,  the representations for the {\sl fundamental solutions of the partial differential operator}. In fact,  this diagram
brings into a single  hierarchy several different partial differential operators.
Indeed, if we take into account the propagation cone by introducing the distance function
$ \phi (t)$, and if we provide   the integral operator (\ref{main}) with the kernel $ K (t;r,b;M) $, as in (\ref{KJDE}),
then we actually  generate  new representations for the solutions of  different well-known equations with $x$-independent coefficients.
 (See, for details, \cite{yagdjian_Rend_Trieste}.)
\smallskip

In the present paper we extend the class of the equations for which we can obtain  explicit representation formulas for the solutions,
by varying the first mapping. 
More precisely, consider the  diagram of Figure~2,
where $w=w_{A,\varphi}(x,t;b)$ is a solution to the  problem (\ref{0.1JDE})
with the parameter $b \in I \subseteq {\mathbb R}$.
If we have a resolving operator of the problem (\ref{0.1JDE}), then,  by applying (\ref{KJDE}),  we can generate   solutions of  another   equation.
Thus, ${\mathcal G}_A = {\mathcal K}\circ {\mathcal E}{\mathcal E}_A $. The new class of equations contains operators with $x$-depending coefficients, and those equations
are not necessarily hyperbolic.

\smallskip

In this paper we restrict ourselves to the  Klein-Gordon  equation in the de~Sitter spacetime,
that is $a(t)=e^{-t}$ in (\ref{0.1JDE}). Recently the  equations in the de Sitter and anti-de Sitter spacetimes
became the focus of interest for an increasing number of authors
(see, e.g., \cite{Bachelot}-\cite{BaskinSE}, \cite{Bros}-\cite{Costa}, \cite{Dohse}, \cite{Hintz}-\cite{Kong-Wei}, \cite{Moschella}
-\cite{Nakamura}, \cite{Rendall_2004,Vasy_2010,Vasy_2013,Helsinki_2013} and the bibliography therein) which investigate those equations from a wide spectrum of perspectives.
 The creation
of a tool for the investigation of the local and global solvability in the  problems for these linear and nonlinear equations appears to
be a worthwhile undertaking. We believe that the integral transform and the representation formulas for the solutions that we derive in this article fill up the gap in the literature on that topic.
\smallskip

To formulate the main result of this paper we need the following notations.
First, we define a  {\it chronological future}
$D_+ (x_0,t_0) $
and  a
{\it chronological past} $D_- (x_0,t_0) $  of the point   $(x_0 ,t_0)$,  $x_0 \in {\mathbb R}^n$, $t_0 \in {\mathbb R}$,
 as follows:
$
D_\pm (x_0,t_0)
  :=
 \{ (x,t)  \in {\mathbb R}^{n+1}  \, ; \,
|x -x_0 | \leq $ $\pm( e^{-t_0} - e^{-t })
\, \} $.
Then, for $(x_0, t_0) \in {\mathbb R}^n\times {\mathbb R}$, $M \in {\mathbb C}$,  we define the function
\begin{eqnarray}
\label{EM}
E(x,t;x_0,t_0;M)
& :=  &
 4 ^{-M}  e^{ M(t_0+t) } \Big((e^{-t _0}+e^{-t})^2 - (x - x_0)^2\Big)^{M -\frac{1}{2}} \\
 &  &
\times F\Big(\frac{1}{2}-M   ,\frac{1}{2}-M  ;1;
\frac{ ( e^{-t_0}-e^{-t })^2 -(x- x_0 )^2 }{( e^{-t_0}+e^{-t })^2 -(x- x_0 )^2 } \Big) , \nonumber
\end{eqnarray}
where   $(x,t) \in D_+ (x_0,t_0)\cup D_- (x_0,t_0) $  and $F\big(a, b;c; \zeta \big) $ is the hypergeometric function. (For definition of the hypergeometric function, see, e.g., \cite{B-E}.)
When no ambiguity arises, like in (\ref{EM}), we use the notation $x^2:= |x|^2$ for $x \in {\mathbb R}^n $.
Thus, the function $E  $ depends on $r^2= (x- x_0 )^2$, that is $E(x,t;x_0,t_0;M)= E  (r,t;0,t_0;M) $. According to Theorem~\ref{T2.12}, the function $  E  (r,t;0,t_0;M) $ solves  the following one dimensional Klein-Gordon equation in the de~Sitter spacetime:
 \begin{eqnarray*}
&  &
E_{tt} (r,t;0,t_0;M) - e^{-2t} E_{rr}(r,t;0,t_0;M)-M^2 E(r,t;0,t_0;M) =0\,.
\end{eqnarray*}

The kernels   $K_0(z,t;M) $   and  $K_0(z,t;M) $ are defined by
\begin{eqnarray}
\label{K0M}
K_0(z,t;M)
&  := &
4 ^ {-M}  e^{ t M}\big((1+e^{-t })^2 - z^2\big)^{  M    } \frac{1}{ [(1-e^{ -t} )^2 -  z^2]\sqrt{(1+e^{-t } )^2 - z^2} } \\
&   &
\times  \Bigg[  \big(  e^{-t} -1 +M(e^{ -2t} -      1 -  z^2) \big)
F \Big(\frac{1}{2}-M   ,\frac{1}{2}-M  ;1; \frac{ ( 1-e^{-t })^2 -z^2 }{( 1+e^{-t })^2 -z^2 }\Big)\nonumber  \\
&  &
\hspace{0.3cm}  +   \big( 1-e^{-2 t}+  z^2 \big)\Big( \frac{1}{2}+M\Big)
F \Big(-\frac{1}{2}-M   ,\frac{1}{2}-M  ;1; \frac{ ( 1-e^{-t })^2 -z^2 }{( 1+e^{-t })^2 -z^2 }\Big) \Bigg]\nonumber ,\\
\label{K1M}
K_1(z,t;M)
& :=  &
  4 ^{-M} e^{ Mt }  \big((1+e^{-t })^2 -   z  ^2\big)^{-\frac{1}{2}+M    }  
  F\left(\frac{1}{2}-M   ,\frac{1}{2}-M  ;1;
\frac{ ( 1-e^{-t })^2 -z^2 }{( 1+e^{-t })^2 -z^2 } \right)\,.  
\end{eqnarray}

The equation (\ref{0.3JDE}) is said to be an equation with imaginary (real) mass if there is  minus (plus) in front of $M^2$; however $M \in {\mathbb C}$.
For the  Klein-Gordon equation with imaginary mass    we have the following result.

\begin{theorem}
\label{T1.2}
For $ f \in C^\infty (\Omega\times I  )$,\, $ I=[0,T]$, $0< T \leq \infty$, and \, $ \varphi_0 $,  $ \varphi_1 \in C_0^\infty (\Omega ) $,
let  the function
$v_f(x,t;b)$
be a solution to the   problem
\begin{equation}
\label{1.22}
\cases{ v_{tt} -   A(x,\partial_x) v  =  0 \,, \quad x \in \Omega \,,\quad t \in [0,1-e^{-T}]\,,\cr
v(x,0;b)=f(x,b)\,, \quad v_t(x,0;b)= 0\,, \quad b \in I,\quad x \in \Omega\,,}
\end{equation}
and the function $v_\varphi = v_\varphi(x,  t)$  be a   solution   of the   problem
\begin{equation}
\label{1.23}
\cases{ v_{tt}-   A(x,\partial_x) v =0, \quad x \in \Omega \,,\quad t \in [0,1-e^{-T}]\,, \cr
 v(x,0)= \varphi (x), \quad v_t(x,0)=0\,,\quad x \in \Omega\,.}
\end{equation}

Then the function  $u= u(x,t)$   defined by
\begin{eqnarray}
\label{1.21}
u(x,t)
&  =  &
2   \int_{ 0}^{t} db
  \int_{ 0}^{\phi (t)- \phi (b)}  v_f(x,r ;b) E(r,t;0,b;M)  \, dr  
+ e ^{\frac{t}{2}} v_{\varphi_0}  (x, \phi (t))\\
&  &
+ \, 2\int_{ 0}^{\phi (t)} v_{\varphi_0}  (x, s) K_0( s,t;M)   ds  \nonumber 
+\, 2\int_{0}^{\phi (t) }  v_{\varphi _1 } (x,  s)
  K_1( s,t;M)   ds
, \quad x \in \Omega  , \,\, t \in I\,,
\end{eqnarray}
where $\phi (t):= 1-e^{-t} $,
 solves the problem
\begin{equation}
\label{1.20}
\cases{u_{tt} - e^{-2t}A(x,\partial_x) u - M^2 u= f, \quad  x \in \Omega \,,\,\, t \in I,\cr
  u(x,0)= \varphi_0 (x)\, , \quad u_t(x,0)=\varphi_1 (x),\quad x \in \Omega\,.}
\end{equation}
Here the kernels  $E$, $K_0$ and $K_1$ have been defined in (\ref{EM}), (\ref{K0M}) and (\ref{K1M}), respectively.
 \end{theorem}

We note    that the operator $A(x,\partial_x)$  is of arbitrary order, that is, the equation of (\ref{1.20}) can be an evolution equation, not necessarily hyperbolic.
Then,   the problems in (\ref{1.22}) and (\ref{1.20}) can be a mixed initial-boundary value problem involving the boundary condition.   Indeed, assume that 
$\Omega \subset {\mathbb R}^n $ is domain with smooth boundary $\partial \Omega $, and that  $ \nu =\nu (x)$ is a unit normal vector. Let $\alpha =\alpha (x) $ and $ \beta =\beta (x)$ be continuous functions, $\alpha ,\beta \in  C(\partial \Omega)$.  If $v=v(t,x) $ satisfies the boundary condition
\[
\alpha (x) v(x,t)+\beta (x) \partial_\nu  v(x,t)=0 \quad \mbox{\rm for all} \quad t\in [0,1-e^T],  \quad x \in \partial \Omega \,,
\]
 then the function  $u= u(x,t)$ fulfills the same boundary condition
\[
\alpha (x) u(x,t)+\beta (x) \partial_\nu  u(x,t)=0 \quad \mbox{\rm for all} \quad t\in I,  \quad x \in \partial \Omega \,.
\]
Next, we stress that interval $[0,1-e^{-T}]\subseteq [0,1]$, which appears in (\ref{1.22}), reflects the fact that de~Sitter model possesses the horizon \cite{Hawking}; existence of the horizon
in the de~Sitter model is widely used to define an asymptotically de Sitter space  (see, e.g., \cite{Baskin,Vasy_2010}) and to involve  geometry into the analysis of the operators on the de~Sitter space (see, e.g., \cite{Birrell,Moschella,Parker-Toms,Tolman}).
\smallskip

Although the next statement  is a straightforward consequence of Theorem~\ref{T1.2}, we present it as a separate theorem having in mind the importance of the equation   with real mass. In the next theorem we use the following kernels
\begin{eqnarray}
\label{Ea}
E(x,t;x_0,t_0)
& := &
E(x,t;x_0,t_0;-iM), \\
\label{K0}
K_0(z,t )&  : = & K_0(z,t;-iM)\,,\\
\label{K1}
K_1(z,t) & : = & K_1(z,t;-iM)\,,
\end{eqnarray}
which were introduced in \cite{Yag_Galst_CMP}.

\begin{theorem}
\label{T1.1}
For \, $ f \in C^\infty (\Omega\times I  )$,\, $ I=[0,T]$, $0< T \leq \infty$,  and
  \, $ \varphi_0 $,  $ \varphi_1 \in C_0^\infty (\Omega ) $, let
 the function
$v_f(x,t;b)$
be a solution to the   problem
\[
\cases{ v_{tt} -   A(x,\partial_x) v  =  0 \,,  \quad x \in \Omega \,,\quad t \in [0,1-e^{-T}]\,,\cr
v(x,0;b)=f(x,b)\,, \quad v_t(x,0;b)= 0\,, \quad b \in I,\quad x \in \Omega\,,}
\]
and the function $v_\varphi = v_\varphi(x, t)$  be a   solution   of the   problem
\[
\cases{ v_{tt}-   A(x,\partial_x) v =0, \quad x \in \Omega \,,\quad t \in [0,1-e^{-T}]\,,\cr
 v(x,0)= \varphi (x), \quad v_t(x,0)=0\,,\quad x \in \Omega\,.}
\]

Then the function  $u= u(x,t)$  defined by
\begin{eqnarray*}
u(x,t)
&  =  &
2   \int_{ 0}^{t} db
  \int_{ 0}^{ \phi (t)- \phi (b)}  v_f(x,r ;b) E(r,t;0,b)  \, dr   
+ e ^{\frac{t}{2}} v_{\varphi_0}  (x, \phi (t))\\
&  &
+ \, 2\int_{ 0}^{\phi (t)} v_{\varphi_0}  (x, s) K_0(s,t)\,  ds   
+\, 2\int_{0}^{\phi (t)}   v_{\varphi _1 } (x, s)
  K_1(s,t) \, ds
, \quad x \in \Omega , \,\, t \in I\,,
\end{eqnarray*}
where $\phi (t):= 1-e^{-t} $, solves  the problem
\[
\cases{u_{tt} - e^{-2t}A(x,\partial_x) u + M^2 u= f, \quad x \in \Omega \,,\,\,t \in I,  \cr
  u(x,0)= \varphi_0 (x)\, , \quad u_t(x,0)=\varphi_1 (x),\quad x \in \Omega\,.}
\]
Here the kernels $E$, $K_0$ and $K_1$ have been defined in (\ref{Ea}), (\ref{K0}) and (\ref{K1}), respectively.
\end{theorem}

 The special cases of Theorems~\ref{T1.2},~\ref{T1.1}, when $A(x,\partial_x)=\Delta  $,  one can find in
\cite{Yag_Galst_CMP,yagdjian_DCDS}.  The proof given in those papers is based on the well-known explicit
representation formulas for the wave equation, the Riemann function, the spherical means,  and the Asgeirsson's mean
value theorem. The main outcome, resulting from  the application of all those tools, is the derivation of the final
representation formula and the kernels $E$, $K_0$, and $K_1$.   Having in the hand the integral transform and the final formulas, we suggest here  
straightforward proof by substitution, which works also for the equations with coefficients depending on $x$.
\smallskip

Among possible applications  of the integral transform method are the $L^p-L^q$ estimates, Strichartz estimates, Huygens' principle, global and local existence theorem for semilinear and quasilinear equations.
Below we give  examples of the equations with the variable coefficients those are amenable to the integral transform method.
\medskip

\noindent
{\bf Example 1.}
The metric $g$ in the de~Sitter type 
spacetime, that is, 
$g_{00}=  g^{00}= -  1  $, $g_{0j}= g^{0j}= 0$, $g_{ij}(x,t)=e^{2t}    \delta _{ij} (x) $, $|g(x,t)|=  e^{2nt}   |\det \delta (x)|$, $g^{ij}(x,t)=  e^{-2t} \delta ^{ij} (x)  $,
$i,j=1,2,\ldots,n$, where $\sum_{j=1}^n\delta ^{ij} (x) \delta _{jk} (x)=\delta _{ik} $,   and $\delta _{ij} $ is  Kronecker's delta.
The linear covariant Klein-Gordon equation      in the coordinates  is
\[
\psi _{tt}
-  \frac{e^{-2t}}{\sqrt{|\det \delta ( x)| }} \sum_{i,j=1}^n  \frac{\partial  }{\partial x^i}\left(  \sqrt{|\det \delta ( x)| }  \delta ^{ij} (x)\frac{\partial  }{\partial x^j}  \psi \right)
+ n   \psi_t +  m^2 \psi
 =
   f \,.
\]
Here $m$ is a physical mass of the particle. If we introduce the new unknown function $u=e^{nt/2}\psi  $, then the equation takes the form of the  Klein-Gordon equation with imaginary mass
\[
u_{tt}
-  \frac{e^{-2t}}{\sqrt{|\det \delta ( x)| }} \sum_{i,j=1}^n  \frac{\partial  }{\partial x^i}\left(  \sqrt{|\det \delta ( x)| }  \delta ^{ij} (x)\frac{\partial  }{\partial x^j}  u \right)
- M^2 u
 =    f \,,
\]
where $-M^2=m^2- \frac{n^2}{4}$ is the  square of the so-called curved (or effective) mass. For the last equation we set
\[
A(x,\partial_x)u =\frac{1}{\sqrt{|\det \delta ( x)| }} \sum_{i,j=1}^n  \frac{\partial  }{\partial x^i}\left(  \sqrt{|\det \delta ( x)| }  \delta ^{ij} (x)
\frac{\partial  }{\partial x^j}  u \right)\,.
\]
If $\Omega =\Pi $ is a non-Euclidean space of constant negative curvature  and the equation of the problems (\ref{1.22}) and (\ref{1.23})
is a non-Euclidean wave equation, then the explicit 
representation formulas are known (see, e.g., \cite{Helgason1984,Lax}) and the Huygens' principle is a consequence of those formulas.    
Thus, for a non-Euclidean wave equation, due to Theorem~\ref{T1.2}, the functions $v_f(x,t;b) $ and $v_\varphi (x,t) $ have  explicit 
representations, and the arguments of \cite{Yag_Galst_CMP,JMP2013}    allow us to derive for the solution $u(x,t)$ of the problem (\ref{1.20}) in the  de~Sitter type metric with   hyperbolic spatial geometry the explicit 
representation, the $L^p-L^q$ estimates, and to examine   the Huygens' principle. Precise statements will be published in the forthcoming paper. 
\medskip

\noindent
{\bf Example 2.}
The Euler-Bernoulli  beam equation with the variable coefficients
 \[
\psi _{tt}
+ e^{-2t} \sum_{i,j=1}^n \partial^2_{x_i} \left( a ^{ij} (x)   \partial^2_{ x_j} \psi \right)
 =
   f \,.
\]
Here $A(x,\partial_x)=\sum_{i,j=1}^n   \partial^2_{x_i}a ^{ij} (x) \partial^2_{ x_j} $ and we assume that $\sum_{i,j=1}^na ^{ij} (x)  \xi _i \xi _j\geq 0$.
 \smallskip

This paper is organized as follows.  In Section~\ref{S2} we study the kernel functions $E$, $K_0$, and $K_1$, and prove several basic properties
of those function. Then,
 in Section \ref{S3}, we prove Theorem~\ref{T1.2}.   Applications of  Theorems~\ref{T1.2},\ref{T1.1}  to   
some equations appearing in electrodynamics and cosmology will be done in a forthcoming paper.

\section{The Kernels of the Integral Transforms}
\label{S2}
\setcounter{equation}{0}
\renewcommand{\theequation}{\thesection.\arabic{equation}}

The proof of Theorem~\ref{T1.2} is straightforward; we just substitute the function of  (\ref{1.21})   in the equation of (\ref{1.20}), and   then check the initial conditions.
It is straightforward, but not short; it requires very long and tedious calculations. In order to make the calculations  more transparent, we  reveal in this section several main properties of the kernel functions.

\subsection{The kernel function  $E(r,t; 0,b;M)$}
\label{SS2.1}

In this subsection we collect some important properties of the kernel $E(x,t; x_0,t_0;M) $.
Although the function $E(x,t; x_0,t_0;M) $ (\ref{EM}) is defined for   $x, x_0 \in {\mathbb R}^n$,   we use it    for $x_0=0 $ and $x =r \in {\mathbb R} $, only.
Consider for $ t,b \in {\mathbb R} $, $0 \leq b\leq t $, $M \in {\mathbb C}$, $r \in [0,    e^{-b}-e^{-t }  ] $ the function
\begin{eqnarray*}
E(r,t; 0,b;M)
& :=  &
 4 ^{-M}  e^{ M(b+t) } \Big((e^{-b }+e^{-t})^2 - r^2\Big)^{M -\frac{1}{2}}
 F\Big(\frac{1}{2}-M   ,\frac{1}{2}-M  ;1;
\frac{ ( e^{-b}-e^{-t })^2 -r^2 }{( e^{-b}+e^{-t })^2 -r^2 } \Big) \,.
\end{eqnarray*}
The following notations are helpful to evaluate derivatives of the kernel functions
\begin{eqnarray*}
 \alpha (r,t,b;M)
& : = &
 4 ^{-M}  e^{ M(b+t) } \Big((e^{-b }+e^{-t})^2 - r^2\Big)^{M  }\,,\\
\beta (r,t,b)
& : = &
 \left((e^{-b }+e^{-t})^2 - r^2\right)^{-\frac{1}{2}}\,,\\
\gamma (r,t,b)
& : =  &
\frac{ ( e^{-b}-e^{-t })^2 -r^2 }{( e^{-b}+e^{-t })^2 -r^2 }\, .
\end{eqnarray*}
 We rewrite the function $E(r,t; 0,b;M) $ in terms of these functions as follows
\begin{eqnarray}
\label{2.24}
\label{2.1}
E(r,t; 0,b;M)
& =  &
\alpha (r,t,b;M) \beta (r,t,b)F\Big(\frac{1}{2}-M   ,\frac{1}{2}-M  ;1;\gamma (r,t,b)\Big) \,.
\end{eqnarray}
For the derivatives of the auxiliary functions $\alpha  $, $ \beta $, and $\gamma  $ we have

\begin{lemma}
\label{L2.2new}
The partial derivatives of the functions $\alpha  $, $ \beta $, and $\gamma  $  are as follows
\begin{eqnarray*}
 \alpha_r (r,t,b;M)
&  = &
-2 Mr  \alpha (r,t,b;M) \beta^2 (r,t,b) \,,\\
\beta_r (r,t,b)
&  = &
 r \beta^3 (r,t,b)\,,\\
\gamma_r (r,t,b)
&  =  &
2 r \beta^2 (r,t,b)\left[ \gamma  (r,t,b) - 1\right]\,,\\
 \alpha_{rr} (r,t,b;M)
&  = &
 -2 M  \alpha (r,t,b;M) \beta^2 (r,t,b)\left[1-2M r^2    \beta^2 (r,t,b)+4 r     \beta^2  (r,t,b)\right] \,,\\
\beta_{rr}  (r,t,b)
&  = &
 \beta^3 (r,t,b)+ 3r^2  \beta^5 (r,t,b)\,,\\
\gamma_{rr}  (r,t,b)
& = &
2 \beta^2 (r,t,b)(1+ 4r^2\beta^2 (r,t,b)) \left[ \gamma  (r,t,b) - 1\right]\,.
 \end{eqnarray*}
 \end{lemma}
 \medskip

 \noindent
 {\bf Proof.} We skip the simple proof of formulas for the first-order derivatives and of $\beta_{rr}  (r,t,b)$. 
For the second-order derivatives we have
 \begin{eqnarray*}
 \alpha_{rr} (r,t,b;M)
& = &
 -2 M  \alpha (r,t,b;M) \beta^2 (r,t,b)+4M^2r^2  \alpha (r,t,b;M)  \beta^4 (r,t,b)\\
&  &
-4Mr^2  \alpha (r,t,b;M)   \beta^4  (r,t,b) \\
& = &
 -2 M  \alpha (r,t,b;M) \beta^2 (r,t,b)\left[1-2M r^2    \beta^2 (r,t,b)+4 r     \beta^2  (r,t,b)\right]\,, \\
\gamma_{rr}  (r,t,b)
& = &
 2  \beta^2 (r,t,b)\left[ \gamma  (r,t,b) - 1\right] +  4 r \beta_r  (r,t,b)\beta  (r,t,b)\left[ \gamma  (r,t,b) - 1\right] \\
&  &
+ 2 r \beta^2 (r,t,b) \gamma_r  (r,t,b)  \\
& = &
 2  \beta^2 (r,t,b)\left[ \gamma  (r,t,b) - 1\right] +  4 r^2  \beta^4  (r,t,b)\left[ \gamma  (r,t,b) - 1\right] \\
&  &
+ 4 r^2  \beta^4  (r,t,b)\left[ \gamma  (r,t,b) - 1\right]  \\
& = &
2 \beta^2 (r,t,b)(1+ 4r^2\beta^2 (r,t,b)) \left[ \gamma  (r,t,b) - 1\right]
.
 \end{eqnarray*}
 Lemma is proven. \hfill $\square$
 \begin{corollary}
\label{C1.2}
 The derivatives, explicitly written,  are as follows:
\begin{eqnarray*}
 \alpha_r (r,t,b;M)
&   = &
-2^{1-2 M} M r e^{M (b+t)} \left(\left(e^{-b}+e^{-t}\right)^2-r^2\right)^{M-1} , \\
\beta_r (r,t,b)
&  = &
  r \left(\left(e^{-b}+e^{-t}\right)^2-r^2\right)^{-3/2}   \,,\\
\gamma_r (r,t,b)
&  =  &
- 8 r e^{3 (b+t)} \left( (e^{b }+e^{ t})^2  -r^2e^{2 (b+t)} \right)^{-2}\,.
 \end{eqnarray*}
 \end{corollary}
 In particular, for $r= e^{-b }-e^{-t}$ we obtain the next values of those functions.
\begin{lemma}
\label{L2.2}
For the functions $\alpha  $, $\beta  $, and $\gamma  $ we have
\begin{eqnarray*}
 \alpha (e^{-b }-e^{-t},t,b;M)
&  = &
1, \quad
\beta (e^{-b }-e^{-t},t,b)
    =  
2^{-1}e^{ \frac{1}{2}(b+t) }, \quad
\gamma (e^{-b }-e^{-t},t,b)
    =   
0,
\end{eqnarray*}
while for their derivatives we have
\begin{eqnarray*}
 \alpha_r (e^{-b }-e^{-t},t,b;M)
&   = &
-  2^{-1}M(e^{-b }-e^{-t})   e^{ b+t } , \\
\beta_r (e^{-b }-e^{-t},t,b)
&   = &
2^{-3} (e^{-b }-e^{-t}) e^{ \frac{3}{2}(b+t) } ,\\
\gamma_r (e^{-b }-e^{-t},t,b)
&   =  &
-2^{-1} (e^{-b }-e^{-t})   e^{  b+t}
.
 \end{eqnarray*}
 \end{lemma}
 \medskip

 \noindent
 {\bf Proof.} It is a simple consequence of the definitions of these functions and the previous lemma.    Lemma is proven. \hfill $\square$
\smallskip

\noindent
Now we turn to the derivatives of the function $  E(r,t; 0,b;M) $. 
\begin{proposition}
\label{P2.3}
The derivative $  E_{r} (r,t; 0,b;M) := \partial_r E(r,t; 0,b;M) $ is given as follows
\begin{eqnarray*}
\partial_r E(r,t; 0,b;M)
& =  &
2\left(\frac{1}{2}-M\right)r  \alpha (r,t,b;M) \beta^3 (r,t,b)  \Bigg[  F \Big(\frac{1}{2}-M   ,\frac{1}{2}-M  ;1;\gamma (r,t,b)\Big) \\
&  &
\hspace{1cm} + (\gamma  (r,t,b) - 1) \left(\frac{1}{2}-M\right)  F\left(\frac{3}{2}-M,\frac{3}{2}-M;2;\gamma (r,t,b)\right) \Bigg]\,.
\end{eqnarray*}
Moreover,
\begin{eqnarray*}
 E_r(0,t; 0,b;M)
& = &
0,\\
 E_r(e^{-b }-e^{-t},t; 0,b;M)
& =  &
2^{-2}  \left(\frac{1}{4}-M^2\right)(e^{t }-e^{b})  e^{ \frac{1}{2}(b+t) } \,.
 \end{eqnarray*}
The second-order derivative $  E_{rr} (r,t; 0,b;M) := \partial^2_r E(r,t; 0,b;M) $ is given by
\begin{eqnarray*}
\partial_r^2 E(r,t; 0,b;M)
& =  &
2\left(\frac{1}{2}-M\right)  \alpha (r,t,b;M) \beta^3 (r,t,b)\left[ 1
+ 2\left(\frac{3}{2} -  M \right)r^2  \beta^2 (r,t,b)\right ] \\
&  &
\times \Bigg[  F \Big(\frac{1}{2}-M   ,\frac{1}{2}-M  ;1;\gamma (r,t,b)\Big) \\
&  &
\hspace{0.5cm} + (\gamma  (r,t,b) - 1) \left(\frac{1}{2}-M\right)  F\left(\frac{3}{2}-M,\frac{3}{2}-M;2;\gamma (r,t,b)\right) \Bigg] \\
&  &
+ \,4\left(\frac{1}{2}-M\right)^2 \left(\frac{3}{2}-M\right)r^2  \alpha (r,t,b;M) \beta^5 (r,t,b)\left( \gamma  (r,t,b) - 1\right) \\
&  &
\times  \Bigg[    F   \Big(\frac{3}{2}-M   ,\frac{3}{2}-M  ;2;\gamma (r,t,b)\Big)
 \\
&  &
\hspace{0.5cm} + (\gamma  (r,t,b) - 1)    \left(\frac{3}{2}-M\right)  \frac{1}{2} F  \left(\frac{5}{2}-M,\frac{5}{2}-M;3;\gamma (r,t,b)\right)\Bigg]\,.
\end{eqnarray*}

\end{proposition}
\medskip

 \noindent
 {\bf Proof.} Indeed, due to \cite[(7) Sec.~2.1.2]{B-E}, we have
\begin{eqnarray*}
&  &
\partial_r E(r,t; 0,b;M) \\
& =  &
\left[ \alpha_r (r,t,b;M) \beta (r,t,b) +\alpha (r,t,b;M) \beta_r (r,t,b)\right]F \Big(\frac{1}{2}-M   ,\frac{1}{2}-M  ;1;\gamma (r,t,b)\Big)\\
&  &
+\alpha (r,t,b;M) \beta (r,t,b) \gamma_r (r,t,b)\left(\frac{1}{2}-M\right)^2 F\left(\frac{3}{2}-M,\frac{3}{2}-M;2;\gamma (r,t,b)\right)\\
& =  &
\Bigg[  -2 Mr  \alpha (r,t,b;M) \beta^3 (r,t,b)  \\
&  &
+\alpha (r,t,b;M)    r \beta^3 (r,t,b) \Bigg]F \Big(\frac{1}{2}-M   ,\frac{1}{2}-M  ;1;\gamma (r,t,b)\Big)\\
&  &
+2 r \alpha (r,t,b;M)   \beta^3 (r,t,b)\left[ \gamma  (r,t,b) - 1\right] \left(\frac{1}{2}-M\right)^2 F\left(\frac{3}{2}-M,\frac{3}{2}-M;2;\gamma (r,t,b)\right) \,.
\end{eqnarray*}
This proves the first formula of the lemma and the statement   for $E_r(e^{-b }-e^{-t},t; 0,b;M)$.
Then, for the second-order derivative $\partial_r^2 E(r,t; 0,b;M)$ we have
\begin{eqnarray*}
&  &
\partial_r^2 E(r,t; 0,b;M) \\
& =  &
2\left(\frac{1}{2}-M\right)  \Big \{ \alpha (r,t,b;M) \beta^3 (r,t,b)+ r \alpha_r (r,t,b;M) \beta^3 (r,t,b)\\
&  &
+ 3r \alpha  (r,t,b;M) \beta_r (r,t,b)\beta^2 (r,t,b)\Big \} \\
&  &
\times \Bigg[  F \Big(\frac{1}{2}-M   ,\frac{1}{2}-M  ;1;\gamma (r,t,b)\Big) \\
&  &
+ (\gamma  (r,t,b) - 1) \left(\frac{1}{2}-M\right)  F\left(\frac{3}{2}-M,\frac{3}{2}-M;2;\gamma (r,t,b)\right) \Bigg] \\
&  &
+ 4\left(\frac{1}{2}-M\right)r^2  \alpha (r,t,b;M) \beta^5 (r,t,b)\left[ \gamma  (r,t,b) - 1\right]\\
&  &
\times  \Bigg[   F_z  \Big(\frac{1}{2}-M   ,\frac{1}{2}-M  ;1;\gamma (r,t,b)\Big)
+   \left(\frac{1}{2}-M\right)  F\left(\frac{3}{2}-M,\frac{3}{2}-M;2;\gamma (r,t,b)\right) \\
&  &
+ (\gamma  (r,t,b) - 1) \left(\frac{1}{2}-M\right)  F_z \left(\frac{3}{2}-M,\frac{3}{2}-M;2;\gamma (r,t,b)\right)\Bigg]  \\
& =  &
2\left(\frac{1}{2}-M\right)  \Big \{ \alpha (r,t,b;M) \beta^3 (r,t,b) -2 Mr^2  \alpha (r,t,b;M) \beta^5 (r,t,b)
+ 3r^2 \alpha  (r,t,b;M)  \beta^5 (r,t,b)\Big \} \\
&  &
\times \Bigg[  F \Big(\frac{1}{2}-M   ,\frac{1}{2}-M  ;1;\gamma (r,t,b)\Big) \\
&  &
+ (\gamma  (r,t,b) - 1) \left(\frac{1}{2}-M\right)  F\left(\frac{3}{2}-M,\frac{3}{2}-M;2;\gamma (r,t,b)\right) \Bigg] \\
&  &
+ 4\left(\frac{1}{2}-M\right)^2 \left(\frac{3}{2}-M\right)r^2  \alpha (r,t,b;M) \beta^5 (r,t,b)\left( \gamma  (r,t,b) - 1\right) \\
&  &
\times  \Bigg[    F   \Big(\frac{3}{2}-M   ,\frac{3}{2}-M  ;2;\gamma (r,t,b)\Big)
 \\
&  &
+ (\gamma  (r,t,b) - 1)    \left(\frac{3}{2}-M\right)  \frac{1}{2} F  \left(\frac{5}{2}-M,\frac{5}{2}-M;3;\gamma (r,t,b)\right)\Bigg]  \,.
\end{eqnarray*}
It is easily seen that the last expression coincides with one given in the statement of the proposition.
This completes the proof of the proposition. \hfill $\square$

In fact,   Corollary~\ref{C1.2} allows us to write $\partial_r E(r,t; 0,b;M)$ and $\partial_r^2 E(r,t; 0,b;M)$ in the explicit form. 
\begin{corollary}
\label{C2.4}
We have
\begin{eqnarray*}
\partial_r E(r,t; 0,b;M)
& = &
2\left(\frac{1}{2}-M\right)r   4 ^{-M}  e^{ M(b+t) } \Big((e^{-t }+e^{-b})^2 - r^2\Big)^{M-\frac{3}{2}  }    \\
&  &
\times  \Bigg[  F \Big(\frac{1}{2}-M   ,\frac{1}{2}-M  ;1;\frac{ ( e^{-b}-e^{-t })^2 -r^2 }{( e^{-b}+e^{-t })^2 -r^2 }\Big) \\
&  &
\hspace{0.5cm}  - \frac{  4 e^{-b-t} }{( e^{-b}+e^{-t })^2 -r^2 }  \left(\frac{1}{2}-M\right)  F\left(\frac{3}{2}-M,\frac{3}{2}-M;2;\frac{ ( e^{-b}-e^{-t })^2 -r^2 }{( e^{-b}+e^{-t })^2 -r^2 }\right) \Bigg]\,,\\
\partial_r^2 E(r,t; 0,b;M)
& =  &
\widetilde{A}(r,t,b;M) F \Big(\frac{1}{2}-M,\frac{1}{2}-M;1;\gamma (r,t,b) \Big)\\
&  &
+\widetilde{B}(r,t,b;M) F \Big(\frac{3}{2}-M,\frac{3}{2}-M;2;\gamma (r,t,b) \Big)\\
&  &
+\widetilde{C}(r,t,b;M) F \Big(\frac{5}{2}-M,\frac{5}{2}-M;3;\gamma (r,t,b) \Big)\,,
\end{eqnarray*}
where
\begin{eqnarray*}
\widetilde{A}(r,t,b;M)
& =  &
2^{1-2 M}  \left(\frac{1}{2}-M\right)e^{M (b+t)}  \left(\left(e^{-b}+e^{-t}\right)^2-r^2\right)^{M-\frac{5}{2}}\left[ \left(e^{-b}+e^{-t}\right)^2 +  \left(2-2M \right)  r^2 \right] \,,\\
\widetilde{B}(r,t,b;M)
& =  &
-2^{3-2 M}\left(\frac{1}{2}-M\right)^2 e^{(M-1)  (b+t)}\left(\left(e^{-b}+e^{-t}\right)^2-r^2\right)^{M-\frac{7}{2}} \left[
 \left(e^{-b}+e^{-t}\right)^2
+  \left(5-4M\right)      r^2  \right] ,\\
\widetilde{C}(r,t,b;M)
& =  &
2^{5-2 M} \left(\frac{1}{2}-M\right)^2 \left(\frac{3}{2}-M\right)^2 e^{(M -2)(b+t)}  r^2       \left(\left(e^{-b}+e^{-t}\right)^2   -r^2\right)^{M-\frac{9}{2}}\,.\\
\end{eqnarray*}
\end{corollary}

\smallskip

Next, we turn to the derivatives with respect to time. We skip the proof of the following lemma, which gives the derivatives of the auxiliary functions $ \alpha $, $ \beta $, and $\gamma  $.

\begin{lemma}
\label{L2.6}
For the first-order derivatives of the functions $ \alpha $, $ \beta $, and $\gamma  $ we have
\begin{eqnarray*}
 \alpha_t (r,t,b;M)
&   = &
M \alpha  (r,t,b;M)  -2  Me^{ -t} \left(e^{-b}+e^{-t}\right) \alpha  (r,t,b;M)\beta^2 (r,t,b)\,, \\
\beta_t (r,t,b)
&   = &
 e^{-t} \left(e^{-b}+e^{-t}\right)\beta^3 (r,t,b) \,, \\
\gamma_t (r,t,b)
&   =  &
  2e^{-t } \beta ^2  (r,t,b) \Big\{ ( e^{-b}-e^{-t })     +     \left(e^{-b}+e^{-t}\right)\gamma  (r,t,b)  \Big\} \,,
\end{eqnarray*}
and
\begin{eqnarray*}
 \alpha_t (e^{-b}-e^{-t },t,b;M)
&   = &
M    -  2^{-1}  M \left(e^{-b}+e^{-t}\right) e^{b }\,, \\
\beta_t (e^{-b}-e^{-t },t,b)
&   = &
2^{-3}e^{-t} \left(e^{-b}+e^{-t}\right)e^{ \frac{3}{2}(b+t) }\,, \\
\gamma_t (e^{-b}-e^{-t },t,b)
&   =  &
  2^{-1} e^{-t }e^{ (b+t) }   ( e^{-b}-e^{-t }) \,.
\end{eqnarray*}
\end{lemma}
\medskip

\begin{proposition}
\label{P2.7}
One can write
\begin{eqnarray*}
\partial_t E(r,t; 0,b;M)
& =  &
A(r,t,b;M)  F \left(\frac{1}{2}-M,\frac{1}{2}-M;1;\gamma (r,t,b) \right)\\
&  &
+B(r,t,b;M)
F \left(\frac{3}{2}-M,\frac{3}{2}-M;2;\gamma (r,t,b)\right)\,,
\end{eqnarray*}
where
\begin{eqnarray*}
A(r,t,b;M)
& = &
\alpha  (r,t,b;M) \beta (r,t,b)\left\{   M  +2\left(\frac{1}{2}-  M \right)e^{-t} \left(e^{-b}+e^{-t}\right)     \beta^2 (r,t,b)  \right\}\,. \\
B(r,t,b;M)
& = &
 2e^{-t }  \alpha (r,t,b;M)    \beta ^3  (r,t,b) \left\{ ( e^{-b}-e^{-t })     +     \left(e^{-b}+e^{-t}\right)\gamma  (r,t,b)  \right\} \,.\\
\end{eqnarray*}
\end{proposition}
\medskip

\noindent
{\bf Proof.} Indeed, it is easily seen that
\begin{eqnarray*} 
\partial_t E(r,t; 0,b;M)  
& =  &
\alpha  (r,t,b;M) \beta (r,t,b)\Bigg(  \Big\{ 1  -2  e^{ -t} \left(e^{-b}+e^{-t}\right)  \beta^2 (r,t,b) \Big\} M   \\
&  &
+  e^{-t} \left(e^{-b}+e^{-t}\right)     \beta^2 (r,t,b)\Bigg) F\Big(\frac{1}{2}-M   ,\frac{1}{2}-M  ;1;\gamma (r,t,b)\Big)  \\
&  &
+    2e^{-t }  \alpha (r,t,b;M)    \beta ^3  (r,t,b) \Big\{ ( e^{-b}-e^{-t })     +     \left(e^{-b}+e^{-t}\right)\gamma  (r,t,b)  \Big\} \\
&  &
\times \left( \frac{1}{2}-M\right)^2 F\Big(\frac{3}{2}-M   ,\frac{3}{2}-M  ;2;\gamma (r,t,b)\Big) \\
& =  &
\alpha  (r,t,b;M) \beta (r,t,b)\Bigg(   M  +2\left(\frac{1}{2}-  M \right)e^{-t} \left(e^{-b}+e^{-t}\right)     \beta^2 (r,t,b)  \Bigg)  \\
&  &
\times F\Big(\frac{1}{2}-M   ,\frac{1}{2}-M  ;1;\gamma (r,t,b)\Big)  \\
&  &
+    2e^{-t }  \alpha (r,t,b;M)    \beta ^3  (r,t,b) \Big\{ ( e^{-b}-e^{-t })     +     \left(e^{-b}+e^{-t}\right)\gamma  (r,t,b)  \Big\} \\
&  &
\times \left( \frac{1}{2}-M\right)^2 F\Big(\frac{3}{2}-M   ,\frac{3}{2}-M  ;2;\gamma (r,t,b)\Big)
\end{eqnarray*}
and, consequently, 
\begin{eqnarray*}
&  &
\partial_t E(r,t; 0,b;M) \\
& =  &
-2^{1-2 M} \left(M-\frac{1}{2}\right) \left(e^{-b}+e^{-t}\right) e^{M (b+t)-t} \left(\left(e^{-b}+e^{-t}\right)^2-r^2\right)^{M-\frac{3}{2}} \\
&  &
\times F \left(\frac{1}{2}-M,\frac{1}{2}-M;1;\gamma (r,t,b) \right)+\\
&  &
+ 4^{-M} M e^{M (b+t)} \left(\left(e^{-b}+e^{-t}\right)^2-r^2\right)^{M-\frac{1}{2}}
F \left(\frac{1}{2}-M,\frac{1}{2}-M;1;\gamma (r,t,b) \right)\\
&  &
+4^{-M} \left(\frac{1}{2}-M\right)^2 e^{M (b+t)} \left(\frac{2 e^{-t} \left(e^{-b}-e^{-t}\right)}{\left(e^{-b}+e^{-t}\right)^2-r^2}+\frac{2 e^{-t} \left(e^{-b}+e^{-t}\right) \left(\left(e^{-b}-e^{-t}\right)^2-r^2\right)}{\left(\left(e^{-b}+e^{-t}\right)^2-r^2\right)^2}\right)\\
&  &
\times  \left(\left(e^{-b}+e^{-t}\right)^2-r^2\right)^{M-\frac{1}{2}}
F \left(\frac{3}{2}-M,\frac{3}{2}-M;2;\gamma (r,t,b) \right)\,.
\end{eqnarray*}
Proposition is proven. \hfill $\square$

\begin{corollary}
\label{C2.8}
The coefficients $A(r,t,b;M) $  and $ B(r,t,b;M)$, explicitly written, are as follows
\begin{eqnarray*}
A(r,t,b;M)
& = &
2^{-2M} e^{M (b+t)} \left(\left(e^{-b}+e^{-t}\right)^2-r^2\right)^{M-\frac{3}{2}}\\
&  &
\times \Bigg[ M \left(\left(e^{-b}+e^{-t}\right)^2-r^2\right)-\left(2  M-1\right) \left(e^{-b}+e^{-t}\right) e^{-t}  \Bigg] \,,\\
B(r,t,b;M)
& = &
2^{1-2M} \left(\frac{1}{2}-M\right)^2 e^{M (b+t)}  e^{-t} \left(\left(e^{-b}+e^{-t}\right)^2-r^2\right)^{M-\frac{5}{2}}\\
&  &
\times \left[ \left(e^{-b}-e^{-t}\right)   \left(e^{-b}+e^{-t}\right)^2-\left(e^{-b}-e^{-t}\right)r^2  \right.\\
&  &
\left.+  \left(e^{-b}+e^{-t}\right)  \left(e^{-b}-e^{-t}\right)^2-\left(e^{-b}+e^{-t}\right) r^2    \right] \,.
\end{eqnarray*}
\end{corollary}
In particular,
\begin{eqnarray}
\label{2.2}
E_t(e^{-b}-e^{-t},t; 0,b;M)
& =  &
\frac{1}{16}   e^{\frac{1}{2}( b-t)}  \left(e^b \left(1-4 M^2\right)+\left(4 M^2+3\right) e^t\right)\,.
\end{eqnarray}
Then, the following proposition will be used.
\begin{proposition}
\label{P2.9}
For all $t,b,M$, $b \leq t$, we have
\begin{eqnarray*}
&  &
2 E_r(e^{-b}- e^{-t},t;0,b;M) + 2e^t E_t( e^{-b}- e^{-t},t;0,b;M)- e^t (4e^{-b-t})^{-\frac{1}{2}}=0\,.
\end{eqnarray*}
\end{proposition}

\noindent
{\bf Proof.} According to Proposition~\ref{P2.3} and (\ref{2.2}),  we can write
\begin{eqnarray*}
&  &
2 E_r(e^{-b}- e^{-t},t;0,b;M) + 2e^t E_t( e^{-b}- e^{-t},t;0,b;M)- e^t (4e^{-b-t})^{-\frac{1}{2}}\\
&  = &
 \frac{1}{8} (4 M^2-1) \left(e^b-e^t\right) e^{M (b+t)} \left(e^{-b-t}\right)^{M-\frac{1}{2}}\\
&  &
+  \frac{1}{8} e^{M(b  +  t)} \left(e^{-b-t}\right)^{M-\frac{1}{2}} \left(e^b \left(1-4 M^2\right)+\left(4 M^2-1\right) e^t+ 4e^t \right)- e^t (4e^{-b-t})^{-\frac{1}{2}}=
0\,.
\end{eqnarray*}
Proposition is proven. \hfill $  \square$
\smallskip

\begin{proposition}
\label{P2.10}
For all $r,t,b,M $, $b \leq t$, we have
\begin{eqnarray*}
\partial_t^2 E(r,t; 0,b;M)
& =  &
  A_t(r,t,b;M)  F \left(\frac{1}{2}-M,\frac{1}{2}-M;1;\gamma (r,t,b) \right) \\
&  &
+  C (r,t,b;M)
F \left(\frac{3}{2}-M,\frac{3}{2}-M;2;\gamma (r,t,b)\right) \\
&  &
+D(r,t,b;M)
F \left(\frac{5}{2}-M,\frac{5}{2}-M;3;\gamma (r,t,b)\right)\,,
\end{eqnarray*}
where
\begin{eqnarray*}
C (r,t,b;M)
& =  &
 A (r,t,b;M) \gamma_t (r,t,b)\left( \frac{1}{2}-M\right)^2 +B_t(r,t,b;M) \,, \\
D(r,t,b;M)
& = &
\frac{1}{2}\left( \frac{3}{2}-M\right)^2B(r,t,b;M)  \gamma_t (r,t,b)\,.
\end{eqnarray*}
\end{proposition}
\medskip

\noindent
{\bf Proof.} It can be verified by simple calculations.
Proposition is proven. \hfill $\square$

\begin{corollary}
\label{C2.11}
The functions $A_t (r,t,b;M) $, $ C (r,t,b;M)$, and $D(r,t,b;M) $, explicitly written, are as follows
\begin{eqnarray*}
&  &
 A_t (r,t,b;M)\\
 & = &
4^{-M} M^2 e^{M (b+t)} \left(\left(e^{-b}+e^{-t}\right)^2-r^2\right)^{M-\frac{1}{2}}\\
&  &
-2^{-2 M} (2 M-1) e^{M (b+t)-b-2 t}\left(\left(e^{-b}+e^{-t}\right)^2-r^2\right)^{M-\frac{3}{2}}\Big[   2 e^b M-2 e^b+2 M e^t-e^t  \Big]\\
&  &
+2^{2-2 M} \left(M-\frac{1}{2}\right) \left(M-\frac{3}{2}\right) \left(e^{-b}+e^{-t}\right)^2 e^{M (b+t)-2 t} \left(\left(e^{-b}+e^{-t}\right)^2-r^2\right)^{M-\frac{5}{2}}\,,\\
&  &
  C (r,t,b;M)\\
&  = &
-2^{3-2 M} \left(\frac{1}{2}-M\right)^2 M e^{M (b+t)-3 b-t} \left(\left(e^{-b}+e^{-t}\right)^2-r^2\right)^{M-\frac{5}{2}}\left[-e^{2 b} r^2-e^{2 b-2 t}+1\right] \\
&  &
-2^{2-2M} \left(\frac{1}{2}-M\right)^2 e^{(M-5) (b+t)}\left(\left(e^{-b}+e^{-t}\right)^2-r^2\right)^{M-\frac{7}{2}}\\
&  &
\times \Big[  -4 M r^2 e^{3 b+3 t}-4 M r^2 e^{4b+2t }+4 M e^{2 (b+t)}-4 M e^{3 b+t}+4 M e^{b+3 t}-4 e^{4 b} M\\
&  &
+r^4 e^{4 b+4 t}+4 r^2 e^{3 b+3 t}-2 r^2 e^{2 b+4 t}+8 r^2 e^{4b+2t }-8 e^{2 (b+t)}-4 e^{b+3 t}+3 e^{4 b}+e^{4 t} \Big]\,,\\
&  &
 D(r,t,b;M) \\
& = &
2^{3-2 M } \left(\frac{1}{2}-M\right)^2 \left(\frac{3}{2}-M\right)^2 e^{M (b+t) -2t} \left(\left(e^{-b}+e^{-t}\right)^2-r^2\right)^{M-\frac{9}{2}}\\
&  &
\times \left[e^{-2 b} r^4+2 r^2 e^{-2 b-2 t}-2 e^{-4 b} r^2-2 e^{-4 b-2 t}+e^{-2 b-4 t}+e^{-6 b} \right] \,.
\end{eqnarray*}
\end{corollary}
\medskip

\begin{theorem}
\label{T2.12}
The function $E  (r,t;0,b;M) $ solves the following equation
 \begin{eqnarray*}
&  &
E_{tt} (r,t;0,b;M) - e^{-2t} E_{rr}(r,t;0,b;M)-M^2 E(r,t;0,b;M) =0\,.
\end{eqnarray*}
\end{theorem}
\medskip

\noindent
{\bf Proof.} We have
\begin{eqnarray*}
&  &
E_{tt} (r,t;0,b;M) - e^{-2t} E_{rr}(r,t;0,b;M)-M^2 E(r,t;0,b;M)  \\
& =  &
 I(r,t,b;M) F \left(\frac{1}{2}-M,\frac{1}{2}-M;1;\gamma (r,t,b) \right)
+  J(r,t,b;M)
F \left(\frac{3}{2}-M,\frac{3}{2}-M;2;\gamma (r,t,b)\right) \\
&  &
+Y(r,t,b;M)F \left(\frac{5}{2}-M,\frac{5}{2}-M;3;\gamma (r,t,b)\right)\,,
\end{eqnarray*}
where
\begin{eqnarray*}
 I(r,t,b;M)
& =  &
 A_t(r,t,b;M) - e^{-2t}  \widetilde{A} (r,t,b;M) - M^2 \alpha (r,t,b;M) \beta  (r,t,b)  \,, \\
J(r,t,b;M)
& = &
C (r,t,b;M)  - e^{-2t}\widetilde{B} (r,t,b;M)  \,, \\
Y(r,t,b;M)
& = &
D(r,t,b;M)  - e^{-2t}\widetilde{C} (r,t,b;M) \,.
\end{eqnarray*}
That is,
\begin{eqnarray*}
 I(r,t,b;M)
& =  &
4^{-M} M^2 e^{M (b+t)} \left(\left(e^{-b}+e^{-t}\right)^2-r^2\right)^{M-\frac{1}{2}}\\
&  &
-2^{-2 M} (2 M-1) e^{M (b+t)-b-2 t}\left(\left(e^{-b}+e^{-t}\right)^2-r^2\right)^{M-\frac{3}{2}}\Big[   2 e^b M-2 e^b+2 M e^t-e^t  \Big]\\
&  &
+2^{2-2 M} \left(M-\frac{1}{2}\right) \left(M-\frac{3}{2}\right) \left(e^{-b}+e^{-t}\right)^2 e^{M (b+t)-2 t} \left(\left(e^{-b}+e^{-t}\right)^2-r^2\right)^{M-\frac{5}{2}} \\
&  &
-e^{-2 t}\left(2^{2-2 M} \left(M-\frac{3}{2}\right) \left(M-\frac{1}{2}\right) r^2 e^{M (b+t)} \left(\left(e^{-b}+e^{-t}\right)^2-r^2\right)^{M-\frac{5}{2}} \right.\\
& &
-\left.2^{1-2 M} \left(M-\frac{1}{2}\right) e^{M (b+t)} \left(\left(e^{-b}+e^{-t}\right)^2-r^2\right)^{M-\frac{3}{2}}\right)\,.
\end{eqnarray*}
After simplification we obtain
\begin{eqnarray*}
 I(r,t,b;M)
& =  &
-4^{-M} (1-2 M)^2 e^{(M-1) (b+t)} \left(  \left(e^{-b}+e^{-t}\right)^2-r^2\right)^{M-\frac{3}{2}}\,.
\end{eqnarray*}
Then
\begin{eqnarray*}
J(r,t,b;M)
& =  &
4^{-M} \left(\frac{1}{2}-M\right)^2 e^{M (b+t)} \left(\left(e^{-b}+e^{-t}\right)^2-r^2\right)^{M-\frac{1}{2}}\\
&    &
\times \left(\frac{8 e^{-2 t} \left(e^{-b}+e^{-t}\right)^2 \left(\left(e^{-b}-e^{-t}\right)^2-r^2\right)}{\left(\left(e^{-b}+e^{-t}\right)^2-r^2\right)^3} \right.\\
&   &
+\frac{8 e^{-2 t} \left(e^{-b}-e^{-t}\right) \left(e^{-b}+e^{-t}\right)}{\left(\left(e^{-b}+e^{-t}\right)^2-r^2\right)^2}-\frac{2 e^{-2 t} \left(\left(e^{-b}-e^{-t}\right)^2-r^2\right)}{\left(\left(e^{-b}+e^{-t}\right)^2-r^2\right)^2} \\
&    &
-\frac{2 e^{-t} \left(e^{-b}+e^{-t}\right) \left(\left(e^{-b}-e^{-t}\right)^2-r^2\right)}{\left(\left(e^{-b}+e^{-t}\right)^2-r^2\right)^2}+\frac{2 e^{-2 t}}{\left(e^{-b}+e^{-t}\right)^2-r^2} \\
&    &
-\left.\frac{2 e^{-t} \left(e^{-b}-e^{-t}\right)}{\left(e^{-b}+e^{-t}\right)^2-r^2}\right) \\
&   &
- 2^{2-2 M} \left(\frac{1}{2}-M\right)^2 \left(M-\frac{1}{2}\right) \left(e^{-b}+e^{-t}\right) e^{M (b+t)-t}\\
&   &
\times \left(\left(e^{-b}+e^{-t}\right)^2-r^2\right)^{M-\frac{3}{2}}\\
&   &
\times \left(\frac{2 e^{-t} \left(e^{-b}-e^{-t}\right)}{\left(e^{-b}+e^{-t}\right)^2-r^2}+\frac{2 e^{-t} \left(e^{-b}+e^{-t}\right) \left(\left(e^{-b}-e^{-t}\right)^2-r^2\right)}{\left(\left(e^{-b}+e^{-t}\right)^2-r^2\right)^2}\right) \\
&  &
+2^{1-2 M} \left(\frac{1}{2}-M\right)^2 M e^{M (b+t)} \left(\left(e^{-b}+e^{-t}\right)^2-r^2\right)^{M-\frac{1}{2}}\\
&  &
\times \left(\frac{2 e^{-t} \left(e^{-b}-e^{-t}\right)}{\left(e^{-b}+e^{-t}\right)^2-r^2}+\frac{2 e^{-t} \left(e^{-b}+e^{-t}\right) \left(\left(e^{-b}-e^{-t}\right)^2-r^2\right)}{\left(\left(e^{-b}+e^{-t}\right)^2-r^2\right)^2}\right)\\
&  &
-e^{-2 t}\left(4^{-M} \left(\frac{1}{2}-M\right)^2 e^{M (b+t)} \left(\left(e^{-b}+e^{-t}\right)^2-r^2\right)^{M-\frac{1}{2}}\right.\\
&  &
\times \left(\frac{8 r^2 \left(\left(e^{-b}-e^{-t}\right)^2-r^2\right)}{\left(\left(e^{-b}+e^{-t}\right)^2-r^2\right)^3}-\frac{8 r^2}{\left(\left(e^{-b}+e^{-t}\right)^2-r^2\right)^2}+\right.\\
&  &
+\left.\frac{2 \left(\left(e^{-b}-e^{-t}\right)^2-r^2\right)}{\left(\left(e^{-b}+e^{-t}\right)^2-r^2\right)^2}-\frac{2}{\left(e^{-b}+e^{-t}\right)^2-r^2}\right)\\
&  &
-2^{2-2 M} \left(\frac{1}{2}-M\right)^2 \left(M-\frac{1}{2}\right) r e^{M (b+t)} \left(\left(e^{-b}+e^{-t}\right)^2-r^2\right)^{M-\frac{3}{2}}\\
&  &
\times \left.\left(\frac{2 r \left(\left(e^{-b}-e^{-t}\right)^2-r^2\right)}{\left(\left(e^{-b}+e^{-t}\right)^2-r^2\right)^2}-\frac{2 r}{\left(e^{-b}+e^{-t}\right)^2-r^2}\right)\right)\,.
\end{eqnarray*}
After simplification we obtain
\begin{eqnarray*}
J(r,t,b;M)
& = &
-4^{-M} (1-2 M)^2 e^{(M-3) (b+t)} \left(\left(e^{-b}+e^{-t}\right)^2-r^2\right)^{M-\frac{5}{2}}\\
&  &
\times  \left[2 M r^2 e^{2 (b+t)}+4 M e^{b+t}-2 e^{2 b} M -2 M e^{2 t} -r^2 e^{2 (b+t)} -6 e^{b+t}+e^{2 b}+e^{2 t}\right]\,.
\end{eqnarray*}
Then
\begin{eqnarray*}
Y(r,t,b;M)
& =  &
2^{-2 M-1} \left(\frac{1}{2}-M\right)^2 \left(\frac{3}{2}-M\right)^2 e^{M (b+t)} \left(\left(e^{-b}+e^{-t}\right)^2-r^2\right)^{M-\frac{1}{2}}\\
&  &
\times \left(\frac{2 e^{-t} \left(e^{-b}-e^{-t}\right)}{\left(e^{-b}+e^{-t}\right)^2-r^2}+\frac{2 e^{-t} \left(e^{-b}+e^{-t}\right) \left(\left(e^{-b}-e^{-t}\right)^2-r^2\right)}{\left(\left(e^{-b}+e^{-t}\right)^2-r^2\right)^2}\right)^2-\\
&  &
-2^{-2 M-1} \left(\frac{1}{2}-M\right)^2 \left(\frac{3}{2}-M\right)^2 e^{M (b+t)-2 t} \left(\left(e^{-b}+e^{-t}\right)^2-r^2\right)^{M-\frac{1}{2}}\\
&  &
\times \left(\frac{2 r \left(\left(e^{-b}-e^{-t}\right)^2-r^2\right)}{\left(\left(e^{-b}+e^{-t}\right)^2-r^2\right)^2}-\frac{2 r}{\left(e^{-b}+e^{-t}\right)^2-r^2}\right)^2\,.
\end{eqnarray*}
In fact, we obtain
\begin{eqnarray*} 
Y(r,t,b;M) 
& =  &
2^{3-2 M} \left(\frac{1}{2}-M\right)^2 \left(\frac{3}{2}-M\right)^2 e^{-2 t} e^{(M-2) (b+t)} \left(\left(e^{-b}+e^{-t}\right)^2-r^2\right)^{M-\frac{9}{2}} \\
&  &
\times \left[-2 r^2 e^{2 t-2 b}+e^{2 t-4 b}-2 e^{-2 b}+r^4 e^{2 t}-2 r^2+e^{-2 t}\right]\,.
\end{eqnarray*}
Furthermore,
\begin{eqnarray*}
J(r,t,b;M)F \left(\frac{3}{2}-M,\frac{3}{2}-M;2;\gamma (r,t,b)\right)
& = &
\widetilde{J}(r,t,b;M)  F_{z} \left(\frac{1}{2}-M,\frac{1}{2}-M;1;\gamma (r,t,b)\right)\,,\\
Y(r,t,b;M)F \left(\frac{5}{2}-M,\frac{5}{2}-M;3;\gamma (r,t,b)\right)
& = &
\widetilde{Y}(r,t,b;M)  F_{zz} \left(\frac{1}{2}-M,\frac{1}{2}-M;1;\gamma (r,t,b)\right)\,,
\end{eqnarray*}
where $ F_{z} \left(a,b;1;z\right) := \frac{d}{d z} F  \left(a,b;1;z\right)$   and the following notations have been used
\begin{eqnarray*}
\widetilde{J}(r,t,b;M)
& = &
J(r,t,b;M)\left( \frac{1}{2}-M\right)^{-2} \,, \\
\widetilde{Y}(r,t,b;M)
& = &
Y(r,t,b;M)2\left( \frac{1}{2}-M\right)^{-2} \left( \frac{3}{2}-M\right)^{-2} \,.
\end{eqnarray*}
Hence
\begin{eqnarray*}
&  &
E_{tt} (r,t;0,b;M) - e^{-2t} E_{rr}(r,t;0,b;M)-M^2 E(r,t;0,b;M)  \\
& =  &
 I(r,t,b;M)  F \left(\frac{1}{2}-M,\frac{1}{2}-M;1;\gamma (r,t,b) \right) 
+ \widetilde{J}(r,t,b;M)  F_{z} \left(\frac{1}{2}-M,\frac{1}{2}-M;1;\gamma (r,t,b)\right)\\
&  &
+\widetilde{Y}(r,t,b;M)2  F_{zz} \left(\frac{1}{2}-M,\frac{1}{2}-M;1;\gamma (r,t,b)\right)\,.
\end{eqnarray*}
Next, if we denote
\begin{eqnarray*}
&  &
G(r,t,b;M)
:= 4^{1-M} e^{(M-1) (b+t)} \left( ( e^{-b}+e^{-t})^2 -r^2\right)^{M-\frac{3}{2}} \,,
\end{eqnarray*}
then
\begin{eqnarray*}
&  &
I(r,t,b;M) F \left(\frac{1}{2}-M,\frac{1}{2}-M;1;\gamma (r,t,b) \right) 
+ \widetilde{J}(r,t,b;M)  F_{z} \left(\frac{1}{2}-M,\frac{1}{2}-M;1;\gamma (r,t,b) \right)\\
&  &
+\widetilde{Y}(r,t,b;M) F_{zz} \left(\frac{1}{2}-M,\frac{1}{2}-M;1;\gamma (r,t,b)\right)\\
& = &
G(r,t,b;M)\Bigg\{ \gamma (1-\gamma ) F_{zz}\left(\frac{1}{2}-M,\frac{1}{2}-M;1;\gamma (r,t,b)  \right) \\
&  &
+ \left(1-\left(1+2 \left(\frac{1}{2}-M\right) \right)\gamma\right)F_{z} \left(\frac{1}{2}-M,\frac{1}{2}-M;1;\gamma (r,t,b) \right) \\
&  &
- \left(\frac{1}{2}-M\right)^2 F  \left(\frac{1}{2}-M,\frac{1}{2}-M;1;\gamma (r,t,b) \right)\Bigg\}\,,
\end{eqnarray*}
since 
\begin{eqnarray*}
&  &
G(r,t,b;M)=\frac{I(r,t,b;M)}{- \left(\frac{1}{2}-M\right)^2}
= \frac{\widetilde{ J}(r,t,b;M)}{\left(1-\left(2 \left(\frac{1}{2}-M\right)+1\right)\gamma \right)}= \frac{\widetilde{Y}(r,t,b;M)}{\gamma (1-\gamma )} \,.
\end{eqnarray*}
The hypergeometric function $F  \left(\frac{1}{2}-M,\frac{1}{2}-M;1;z  \right) $ solves the following equation
\begin{eqnarray*}
&  &
z  (1-z  ) F_{zz }
+ \left(1-\left( 2 \left(\frac{1}{2}-M\right)+1 \right)z \right)F_{z }
- \left(\frac{1}{2}-M\right)^2 F  =0\,.
\end{eqnarray*}
 Theorem is proven.
\hfill $\square$

\subsection{The kernel function  $K_1(r,t;M)$}
\label{SS2.2}

By definition
\begin{eqnarray*}
K_1(r,t;M)
& :=  &
  4 ^{-M} e^{ Mt }  \big((1+e^{-t })^2 -   r^2\big)^{-\frac{1}{2}+M    }
  F\left(\frac{1}{2}-M   ,\frac{1}{2}-M  ;1;
\frac{ ( 1-e^{-t })^2 -r^2 }{( 1+e^{-t })^2 -r^2 } \right)\,,
\end{eqnarray*}
where $t \geq 0 $, $r \in [0,1-e^{-t }] $, and $M \in {\mathbb C} $. In fact,
\begin{eqnarray*}
K_1(r,t;M)
&  =  &
 E(r,t;0,0,M)  \,,
\end{eqnarray*}
which simplifies the proof of many properties of this function since they are inherited from the kernel $E(z,t;0,b,M) $.
\begin{proposition}
\label{P2.13}
We have
 \begin{eqnarray*}
&  &
2 e^{-t}     K_{1r}(\phi (t),t;M)  + 2 K_{1t}( \phi (t),t;M)
 -  K_1(\phi (t),t;M) =0 \quad \mbox{\rm for all}  \quad t >0 \,,\\
&  &
K_{1r}( 0,t;M)=0 \quad \mbox{\rm for all}  \quad t >0\,.
\end{eqnarray*}
\end{proposition}
\medskip

\noindent
{\bf Proof.}
Due to  Corollary~\ref{C2.4}, the derivative $\partial_z K_1(z,t;M) $, explicitly written, is as follows
\begin{eqnarray*}
\partial_r K_1(r,t;M)
&  =  &
2\left(\frac{1}{2}-M\right)r  4 ^{-M}  e^{ Mt } \Big((e^{-t }+1)^2 - r^2\Big)^{M  -\frac{3}{2}}     \nonumber  \\
&  &
 \Bigg[  F \Big(\frac{1}{2}-M   ,\frac{1}{2}-M  ;1;\gamma (r,t,0)\Big) \nonumber  \\
&  &
\hspace{1cm} -\frac{  4e^{-t } }{( 1+e^{-t })^2 -r^2 }   \left(\frac{1}{2}-M\right)  F\left(\frac{3}{2}-M,\frac{3}{2}-M;2;\gamma (r,t,0)\right) \Bigg] \nonumber \,.
\end{eqnarray*}
In particular, this proves the second equation of the proposition. Moreover, if $r= \phi (t)=1-e^{-t }$, then
\begin{eqnarray*}
  K_{1\,r} (1-e^{-t },t;M)
&  =  &
\frac{1}{4}\left(\frac{1}{4}-M^2\right)  (1-e^{-t })  e^{\frac{3}{2}t }     \,. \nonumber
\end{eqnarray*}

Then, due to Proposition~\ref{P2.7} and Corollary~\ref{C2.8}, the derivative $ \partial_t K_1(z,t;M)$ is
\begin{eqnarray*}
\partial_t K_1(r,t;M)
& =  &
A(r,t,0;M)  F \left(\frac{1}{2}-M,\frac{1}{2}-M;1;\gamma (r,t,0) \right)\\
&  &
+B(r,t,0;M)
F \left(\frac{3}{2}-M,\frac{3}{2}-M;2;\gamma (r,t,0)\right)\,,
\end{eqnarray*}
where
\begin{eqnarray*}
A (r,t,0;M)
& = &
4^{-M} e^{(M-2) t} \left(\left(e^{-t}+1\right)^2-z^2\right)^{M-\frac{3}{2}} \left(-M  e^{2 t}r^2+M e^{2 t}-M+e^t+1\right)
\end{eqnarray*}
and
\begin{eqnarray*}
B (r,t,0;M)
& = &
-4^{1-M} \left(\frac{1}{2}-M\right)^2 e^{(M-3) t } \left(e^{2 t} \left(r^2-1\right)+1\right) \left(\left(e^{-t}+1\right)^2-r^2\right)^{M-\frac{5}{2}}\,.
\end{eqnarray*}
Hence,
\begin{eqnarray*}
K_{1t} (r,t;M)
& = & 
 4^{-M} e^{(M-2) t} \left(\left(e^{-t}+1\right)^2-r^2\right)^{M-\frac{3}{2}} \left(-M  e^{2 t}r^2+M e^{2 t}-M+e^t+1\right)   \nonumber  \\
&  &
\times F \left(\frac{1}{2}-M,\frac{1}{2}-M;1;\gamma (r,t,0) \right)  \nonumber\\
&  &
  -4^{1-M} \left(\frac{1}{2}-M\right)^2 e^{(M-3) t } \left(e^{2 t} \left(r^2-1\right)+1\right) \left(\left(e^{-t}+1\right)^2-r^2\right)^{M-\frac{5}{2}}  \nonumber\\
&  &
\times
F \left(\frac{3}{2}-M,\frac{3}{2}-M;2;\gamma (r,t,0)\right)\,. \nonumber
\end{eqnarray*}
In particular,
\begin{eqnarray*} 
K_{1\,t} (1-e^{-t},t;M) 
&  =  &
 \frac{1}{16} e^{-\frac{1}{2}t} \left(\left(4 M^2+3\right) e^t-4 M^2+1\right) \,,\nonumber
\end{eqnarray*}
while
\begin{eqnarray*}
K_1(1-e^{-t},t;M)
&  =  &
   \frac{1}{2}      e^{\frac{1}{2} t  } \,.
\end{eqnarray*}
This completes the proof of the proposition. \hfill $\square$
\smallskip

Furthermore, the following proposition easily follows from Theorem~\ref{T2.12}.
\begin{proposition}
\label{P2.14}
The function $ $ solves the following equation
 \begin{eqnarray*}
&  &
  K_{1\,tt}( r,t;M)   - e^{-2t}    K_{1\,rr}(r,t;M)   - M^2
  K_{1 }( r,t;M)  =0 \quad \mbox{\rm for all}  \quad t >0, \quad 0<r<t.
\end{eqnarray*}
\end{proposition}

\subsection{The kernel function $K_{0 }(r,t;M)  $}
\label{SS2.3}

The function $K_0(r,t;M)$ can be written as follows
\begin{eqnarray*}
&  &
K_0(r,t;M) = -\left[ \frac{\partial}{\partial b} E(r,t;0,b;M)\right]_{b=0}\,.
\end{eqnarray*}
We consider $K_0(r,t;M)$ for  $r \in [0, 1-e^{-t })$, and $M \in {\mathbb C}$.
We have
\begin{eqnarray}
\label{1.8JDE}
 K_{0 }(\phi (t),t;M)
& := &
\lim_{r\to\phi (t) }K_{0 }(r,t;M)=
-\frac{1}{4} M^2    e^{t\frac{1}{2}} +\frac{1}{4} M^2 e^{ \frac{3}{2}t} -\frac{3}{16}     e^{ \frac{1}{2}t} -\frac{1}{16}  e^{ \frac{3}{2}t}\,.
\end{eqnarray}
Then
\begin{eqnarray*}
&  &
\partial_r K_0(r,t;M)= K_{0\,r}(r,t;M)\\
& = &
\frac{1}{\left(e^t \left(\left(r^2-1\right) e^t+2\right)-1\right)^2} 4^{-M} (2 M-1) r e^{M t} \left((1+e^{-t})^2-r^2\right)^{M-\frac{3}{2}} \\
&  &
 \times \Bigg\{\Big[M r^4 e^{4 t}+2 M r^2 e^{3 t}-2 M r^2 e^{4 t}+2 M e^t-2 M e^{3 t}+M e^{4 t}-M-3 r^2 e^{2 t}\\
&  &
\hspace{2cm} +r^2 e^{3 t}+e^t-3 e^{2 t}-e^{3 t}+3\Big]    F \left(\frac{1}{2}-M,\frac{1}{2}-M;1;\gamma (r,t,0)\right) \\
&  &
-\frac{2 r^4 e^{5 t}+12 r^2 e^{3 t}-4 r^2 e^{5 t}-14 e^t+12 e^{3 t}+2 e^{5 t} }{e^t \left(\left(r^2-1\right) e^t-2\right)-1}
F \left(\frac{1}{2}-M,\frac{3}{2}-M;1;\gamma (r,t,0)\right)\Bigg\}\,.
\end{eqnarray*}
 Thus,
\begin{eqnarray}
\label{1.9JDE}
 K_{0\,r}(\phi (t),t;M)
  :=  
\lim_{r \to \phi (t)}K_{0\,r}(r,t;M) 
 &  =  &
-\frac{1}{16} M^4   e^{ \frac{1}{2}t}+\frac{1}{8} M^4 e^{ \frac{3}{2}  t}-\frac{1}{16} M^4 e^{ \frac{5}{2}  t}-\frac{7}{32} M^2   e^{ \frac{1}{2}t}  \\
&   &
+\frac{1}{16} M^2 e^{ \frac{3}{2}  t}+\frac{5}{32} M^2 e^{ \frac{5}{2}  t}+\frac{15}{256}    e^{ \frac{1}{2}t}
-\frac{3}{128} e^{ \frac{3}{2}  t} -\frac{9}{256} e^{ \frac{5}{2}  t}   \,.\nonumber
\end{eqnarray} 
Similarly,
\begin{eqnarray*}
&  &
K_{0\,t}(r,t;M)\\
& = &
\frac{4^{-M} e^{(M-2) t} \left(-r^2+e^{-2 t}+2 e^{-t}+1\right)^{M-\frac{3}{2}}}{\left(e^t \left(\left(r^2-1\right) e^t+2\right)-1\right)^2}\nonumber\\
&  &
 \times  \Bigg\{ \Big[M^2 r^6 e^{6 t}-M^2 r^4 e^{4 t}+2 M^2 r^4 e^{5 t}-M^2 r^4 e^{6 t}-M^2 r^2 e^{2 t}+2 M^2 r^2 e^{4 t} \nonumber \\
&  &
 -M^2 r^2 e^{6 t}-2 M^2 e^t-M^2 e^{2 t}+4 M^2 e^{3 t}-M^2 e^{4 t}-2 M^2 e^{5 t}+M^2 e^{6 t}+M^2-M r^4 e^{4 t}\nonumber\\
&  &
-M r^4 e^{6 t}+2 M r^2 e^{2 t}-12 M r^2 e^{4 t}+2 M r^2 e^{6 t}+M e^{2 t}+M e^{4 t}-M e^{6 t}-M-r^2 e^{3 t}\nonumber\\
&  &
  +6 r^2 e^{4 t}-r^2 e^{5 t}+e^t-2 e^{3 t}+e^{5 t}\Big] F\left(\frac{1}{2}-M,\frac{1}{2}-M;1;\frac{\left(1-e^{-t}\right)^2-r^2}{\left(e^{-t}+1\right)^2-r^2}\right)\nonumber\\
&  &
 -\frac{1}{e^t \left(\left(r^2-1\right) e^t-2\right)-1}\Bigg[4 M r^4 e^{5 t}+4 M r^4 e^{7 t}-8 M r^2 e^{3 t}+48 M r^2 e^{5 t}\nonumber\\
&  &
-8 M r^2 e^{7 t}+4 M e^t-4 M e^{3 t}-4 M e^{5 t}+4 M e^{7 t}-2 r^4 e^{5 t}-2 r^4 e^{7 t}+4 r^2 e^{3 t}\nonumber\\
&  &
-24 r^2 e^{5 t}+4 r^2 e^{7 t}-2 e^t+2 e^{3 t}+2 e^{5 t}-2 e^{7 t}\Bigg]  F\left(\frac{1}{2}-M,\frac{3}{2}-M;1;\frac{\left(1-e^{-t}\right)^2-r^2}{\left(e^{-t}+1\right)^2-r^2}\right) \Bigg\}\nonumber
\end{eqnarray*}
and
\begin{eqnarray}
\label{1.10JDE}
 K_{0\,t}(\phi (t),t;M)
& := &
\lim_{r \to \phi (t)}K_{0\,t}(r,t;M)\\
& = &
 \frac{1}{16} M^4   e^{-\frac{1}{2}t} -\frac{1}{8} M^4  e^{\frac{1}{2}t}+\frac{1}{16} M^4 e^{ \frac{3}{2}  t}\nonumber\\
&  &
+\frac{7}{32} M^2   e^{-\frac{1}{2}t} -\frac{3}{16} M^2 e^{\frac{1}{2}t}+\frac{7}{32} M^2 e^{ \frac{3}{2}  t}
+\frac{15}{256}  e^{-\frac{1}{2}t} -\frac{9}{128}   e^{\frac{1}{2}t} -\frac{15}{256} e^{ \frac{3}{2}  t}  \,. \nonumber
\end{eqnarray}
\smallskip

\begin{proposition}
\label{P2.15}
For all $M$ we have
 \begin{eqnarray*}
&  &
K_{0\,r}(0,t;M)=0,\quad \mbox{\rm for all}  \quad t >0, \\
&  &
  K_{0\,tt}( r,t;M)   - e^{-2t}    K_{0\, rr}( r,t;M)   - M^2
  K_{0}( r,t;M)  =0 \quad \mbox{\rm for all}  \quad t >0,\quad r \in[0,1-e^{-t}]\,,\\
&  &
\left( \frac{1}{4}- M^2\right) e ^{\frac{t}{2}}   - 2 e ^{-t}  K_0(\phi (t),t;M) + 4  e ^{-2t} K_{0\,r}(\phi (t),t;M)
+ 4 e ^{-t} K_{0\,t}(\phi (t),t;M)
=0\,.
\end{eqnarray*}
\end{proposition}
\medskip

\noindent
{\bf Proof.} The first and second statements follow from the corresponding ones for the function $E$. Then,  according to (\ref{1.8JDE}), (\ref{1.9JDE}),
and (\ref{1.10JDE}), we have
 \begin{eqnarray*}
&  &
 - 2 e ^{-t}  K_0(\phi (t),t;M) + 4  e ^{-2t} K_{0\,s}(\phi (t),t;M)
+ 4 e ^{-t} K_{0\,t}(\phi (t),t;M)\\
& = &
- 2 e ^{-t} \Bigg[  -\frac{1}{4} M^2    e^{t\frac{1}{2}} +\frac{1}{4} M^2 e^{ \frac{3}{2}t} -\frac{3}{16}     e^{ \frac{1}{2}t} -\frac{1}{16}  e^{ \frac{3}{2}t} \Bigg]    \\
&  &
+ 4  e ^{-2t} \Bigg[  -\frac{1}{16} M^4   e^{ \frac{1}{2}t}+\frac{1}{8} M^4 e^{ \frac{3}{2}  t}-\frac{1}{16} M^4 e^{ \frac{5}{2}  t}-\frac{7}{32} M^2   e^{ \frac{1}{2}t} \\
&   &
+\frac{1}{16} M^2 e^{ \frac{3}{2}  t}+\frac{5}{32} M^2 e^{ \frac{5}{2}  t}+\frac{15}{256}    e^{ \frac{1}{2}t}
-\frac{3}{128} e^{ \frac{3}{2}  t} -\frac{9}{256} e^{ \frac{5}{2}  t}   \Bigg]  \\
&  &
+ 4 e ^{-t} \Bigg[  \frac{1}{16} M^4   e^{-\frac{1}{2}t} -\frac{1}{8} M^4  e^{\frac{1}{2}t}+\frac{1}{16} M^4 e^{ \frac{3}{2}  t}\\
&  &
+\frac{7}{32} M^2   e^{-\frac{1}{2}t} -\frac{3}{16} M^2 e^{\frac{1}{2}t}+\frac{7}{32} M^2 e^{ \frac{3}{2}  t}
-\frac{15}{256}   e^{-\frac{1}{2}t} -\frac{9}{128}   e^{\frac{1}{2}t} -\frac{15}{256} e^{ \frac{3}{2}  t} \Bigg] \\
& = &
- \left( \frac{1}{4}- M^2\right) e ^{\frac{t}{2}}\,.
\end{eqnarray*}
This proves the last statement.
Proposition is proven. \hfill $\square$

\section{Proof of Theorem~\ref{T1.2}}
\label{S3}

\setcounter{equation}{0}
\renewcommand{\theequation}{\thesection.\arabic{equation}}

The proof is straightforward; we just substitute the function $u=u(t,x)$ (\ref{1.21}) into equation (\ref{1.20}) and then check the initial conditions.
In order to make proof more transparent we split it into three independent cases.
\bigskip

\noindent
{\bf Case of ($\varphi_1$) }
In this case set
$
 f(t,x)= 0$, $\,\, \varphi_0 (x) =0
$.
Henceforth   we suppress the  subindex $\varphi_1 $ of $v_{\varphi_1} $. We have
 \begin{eqnarray*}
\frac{1}{2} u(x,t)
& = &
   \int_{0}^{\phi (t)}   v_{ } (x,  s)
   K_1( s,t;M)  \, ds
  , \quad x \in {\mathbb R}^n, \,\, t>0\,,
\end{eqnarray*}
 and consider the derivative
\begin{eqnarray*}
&  &
\frac{1}{2}\partial_t^2 u(x,t) \\
& = &
  - e^{-t}  v_{ } (x,  \phi (t))   K_1(\phi (t),t;M)
+  e^{-2t}  v_{ t} (x,  \phi (t))   K_1(\phi (t),t;M)
+  e^{-2t}  v_{ } (x,  \phi (t))   K_{1s}(\phi (t),t;M) \\
&  &
 +  2 e^{-t}  v_{ } (x,\phi (t))
  K_{1t}( \phi (t),t;M)  +  \int_{0}^{\phi (t)}   v_{ } (x,  s)
  K_{1tt}( s,t;M)  \, ds\,.
\end{eqnarray*}
According to the choice (\ref{1.23}) of the function $v $, we have
\[
\frac{1}{2}A(x,\partial_x) u(x,t)
=
\,   \int_{0}^{\phi (t)}  A(x,\partial_x) v (x,  s)
  K_1( s,t;M)  \, ds  =
\,  \int_{0}^{\phi (t)}     v_{tt } (x,   s)
  K_1( s,t;M)  \, ds   \,.
\]
Now we take into account  (\ref{1.23}), that is $v_{t  } (x, 0)=0 $, and the second statement
$
 K_{1s}( 0,t;M)=0
$
of Proposition~\ref{P2.13}. Hence we have
\[
\frac{1}{2}A(x,\partial_x) u(x,t)
  =  
\,  v_{t  } (x, \phi (t) )
  K_1(\phi (t) ,t;M)
-    v (x,  \phi (t))
  K_{1s}(\phi (t),t;M)    
+  \int_{0}^{\phi (t)}     v  (x,   s)
  K_{1ss}( s,t;M)  \, ds .
\] 
Due to Proposition~\ref{P2.13} and Proposition~\ref{P2.14} we can write
\begin{eqnarray*}
&  &
 \partial_t^2 u(x,t) - e^{-2t} A(x,\partial_x) u(x,t) -M^2u(x,t) \\
& = &
 2e^{-t}  v_{ } (x,  \phi (t))  \Bigg(  -  K_1(\phi (t),t;M)
+  2 e^{-t}     K_{1s}(\phi (t),t;M)  + 2 K_{1t}( \phi (t),t;M)  \Bigg)  \\
  &  &
 +  2\int_{0}^{\phi (t)}   v_{ } (x,  s)  \Bigg(
  K_{1tt}( s,t;M)  \, ds  - e^{-2t}
  K_{1ss}( s,t;M)  \, ds  - M^2
  K_{1 }( s,t;M)   \Bigg) \, ds\,.
\end{eqnarray*}
This competes the proof of the case ($\varphi _1 $).
Similarly, we consider the next case.
\smallskip

\noindent
{\bf Case of ($\varphi_0$) } In this case we set
$
 f(t,x)= 0$, $\, \varphi_1 (x)=0$,
 then (\ref{1.21}) reads
\begin{eqnarray*}
u(x,t)
&  =  &
e ^{\frac{t}{2}} v_{\varphi_0}  (x, \phi (t))
+ \, 2\int_{ 0}^{\phi (t)} v_{\varphi_0}  (x, s) K_0(s,t;M)\,  ds \,.
\end{eqnarray*}
From now on we skip subindex $\varphi_0 $ of $v_{\varphi_0}$.  For the first and second-oder derivatives we have
\begin{eqnarray*}
\partial_t u(x,t)
&  =  &
  \frac{1}{2}e ^{\frac{t}{2}} v_{}  (x, \phi (t)) + e ^{-\frac{t}{2}} v_{t}  (x, \phi (t))\\
 &  &
+ \, 2 e ^{-t}v_{}  (x, \phi (t) ) K_0(\phi (t),t;M) + \, 2\int_{ 0}^{\phi (t)} v_{}  (x, s) K_{0\,t}(s,t;M)\,  ds   \,,
\end{eqnarray*}
and
\begin{eqnarray*}
\partial_t^2 u(x,t)
&  =  &
\frac{1}{4}e ^{\frac{t}{2}} v_{}  (x, \phi (t))
+ e ^{-\frac{3t}{2}} v_{tt}  (x, \phi (t))\\
 &  &
 \, - 2 e ^{-t}v_{}  (x, \phi (t) ) K_0(\phi (t),t;M) + 2 e ^{-\frac{3t}{2}}v_{t}  (x, \phi (t) ) K_0(\phi (t),t;M)\\
&  &
+ 2 e ^{-\frac{3t}{2}}v_{}  (x, \phi (t) ) K_{0\,s}(\phi (t),t;M)+ 2 e ^{-t}v_{}  (x, \phi (t) ) K_{0\,t}(\phi (t),t;M)\\
&  &
+ \, 2 e ^{-t} v_{}  (x,  \phi (t)) K_{0\,t}(\phi (t),t;M)  + \, 2\int_{ 0}^{\phi (t)} v_{}  (x, s) K_{0\,tt}(s,t;M)\,  ds\,,
\end{eqnarray*}
respectively, while
\begin{eqnarray*}
A(x,\partial_x) u(x,t)
&  =  &
e ^{\frac{t}{2}} A(x,\partial_x) v_{ }  (x, \phi (t))
+ \, 2\int_{ 0}^{\phi (t)}   v_{tt }  (x, s) K_0(s,t;M)\,  ds\,.
\end{eqnarray*}
Consider the last integral; using   twice integration by parts and the choice of function $ v$, we obtain
\begin{eqnarray*}
\int_{ 0}^{\phi (t)}   v_{ ss}  (x, s) K_0(s,t;M)\,  ds
&  =  &
  v_{ s }  (x, \phi (t)) K_0(\phi (t),t;M) -     v_{   }  (x, \phi (t)) K_{0\,s}(\phi (t),t;M) \\
&  &
+ v_{   }  (x, 0) K_{0\,s}(0,t;M)+ \int_{ 0}^{\phi (t)}   v_{   }  (x, s) K_{0\,ss}(s,t;M)\,  ds\,.
\end{eqnarray*}
Since $K_{0\,s}(0,t;M)=0 $, we obtain
\begin{eqnarray*}
A(x,\partial_x) u(x,t)
&  =  &
e ^{\frac{t}{2}}   v_{tt }  (x, \phi (t))+    2   v_{ t }  (x, \phi (t)) K_0(\phi (t),t;M) -  2    v_{   }  (x, \phi (t)) K_{0\,s}(\phi (t),t;M) \\
&  &
+2\int_{ 0}^{\phi (t)}   v_{   }  (x, s) K_{0\,ss}(s,t;M)\,  ds\,.
\end{eqnarray*}
Hence
\begin{eqnarray*}
&  &
\partial_t^2 u(x,t) - e ^{-2t}A(x,\partial_x) u(x,t) -M^2 u(x,t)\\
&  =  &
v   (x, \phi (t)) \Bigg[ \left( \frac{1}{4}- M^2\right) e ^{\frac{t}{2}}   - 2 e ^{-t}  K_0(\phi (t),t;M) + 2 e ^{-2t} K_{0\,s}(\phi (t),t;M)\\
&  &
+ 2 e ^{-t}v_{}  (x, \phi (t) ) K_{0\,t}(\phi (t),t;M) 
+ \, 2 e ^{-t}  K_{0\,t}(\phi (t),t;M)
+   2e ^{-2t}     K_{0\,s}(\phi (t),t;M) \Bigg]   \\
  &  &
+2\int_{ 0}^{\phi (t)} v_{}  (x, s) \Bigg[ K_{0\,tt}(s,t;M)\,  ds
-e ^{-2t}  K_{0\,ss}(s,t;M)\,  ds  -   M^2   K_0(s,t;M)\Bigg] \,  ds\,.
\end{eqnarray*}
Then we use Proposition~\ref{P2.15}, which
  completes the  proof of the case of $(\varphi _0) $.
\bigskip

\noindent
{\bf Case of (f)} Similarly, we set
$
\varphi_0 (x)=0 $,  $\,\varphi_1 (x)=0$,
 then (\ref{1.21}) reads
\begin{eqnarray*}
u(x,t)
&  =  &
2   \int_{ 0}^{t} db
  \int_{ 0}^{ e^{-b}- e^{-t}}  v(x,r ;b) E(r,t;0,b;M)  \, dr\,.
\end{eqnarray*}
For the first and second-order derivatives we obtain
\begin{eqnarray*}
\frac{1}{2} \partial_t u(x,t)
&  =  &
 e^{-t}    \int_{ 0}^{t}
v(x, e^{-b}- e^{-t} ;b) E( e^{-b}- e^{-t},t;0,b;M)  \, db  \\
&  &
+  \int_{ 0}^{t} db
  \int_{ 0}^{ e^{-b}- e^{-t}}  v(x,r ;b) E_t(r,t;0,b;M)  \, dr
\end{eqnarray*}
and
\begin{eqnarray*}
\frac{1}{2} \partial_t^2 u(x,t)
&  =  &
-e^{-t}    \int_{ 0}^{t}
v(x, e^{-b}- e^{-t} ;b) E( e^{-b}- e^{-t},t;0,b;M)  \, db   \\
&  &
+e^{-t}
v(x, 0 ;t) E( 0,t;0,t;M)    \\
&  &
+ e^{-2t}    \int_{ 0}^{t}
 v_t(x, e^{-b}- e^{-t} ;b) E( e^{-b}- e^{-t},t;0,b;M)  \, db  \\
&  &
+ e^{-2t}       \int_{ 0}^{t}
v(x, e^{-b}- e^{-t} ;b) E_r( e^{-b}- e^{-t},t;0,b;M)  \, db   \\
&  &
+2e^{-t}    \int_{ 0}^{t}
v(x, e^{-b}- e^{-t} ;b) E_t( e^{-b}- e^{-t},t;0,b;M)  \, db    \\
&  &
+  \int_{ 0}^{t} db
  \int_{ 0}^{ e^{-b}- e^{-t}}  v(x,r ;b) E_{tt} (r,t;0,b;M)  \, dr
\end{eqnarray*}
respectively.
On the other hand, due to the definition of the kernel $ E$,  Lemma~\ref{L2.2}, and (\ref{2.1}), we obtain
\begin{eqnarray}
\label{3.10}
 E( e^{-b}- e^{-t},t;0,b;M)=
\frac{1}{2} e^{\frac{1}{2}(b+t)}  \,, \qquad E( 0,t;0,t;M)=  \frac{1}{2}e^{t}\,.
\end{eqnarray}
Hence, since $ v(x, 0 ;t)=f(x,t) $, we have
\begin{eqnarray*}
\frac{1}{2} \partial_t^2 u(x,t)
&  =  &
-e^{-t}    \int_{ 0}^{t}
v(x, e^{-b}- e^{-t} ;b) (4e^{-b-t})^{-\frac{1}{2}}  \, db
+
\frac{1}{2}  f(x, t)   \\
&  &
+ e^{-2t}    \int_{ 0}^{t}
 v_t(x, e^{-b}- e^{-t} ;b)(4e^{-b-t})^{-\frac{1}{2}}  \, db  \\
&  &
+ e^{-2t}       \int_{ 0}^{t}
v(x, e^{-b}- e^{-t} ;b) E_r( e^{-b}- e^{-t},t;0,b;M)  \, db   \\
&  &
+2e^{-t}    \int_{ 0}^{t}
v(x, e^{-b}- e^{-t} ;b) E_t( e^{-b}- e^{-t},t;0,b;M)  \, db    \\
&  &
+  \int_{ 0}^{t} db
  \int_{ 0}^{ e^{-b}- e^{-t}}  v(x,r ;b) E_{tt} (r,t;0,b;M)  \, dr\,.
\end{eqnarray*}
Then, taking into account (\ref{1.22}), we have
\begin{eqnarray*}
\frac{1}{2} A(x,\partial_x) u(x,t)
&  =  &
 \int_{ 0}^{t} db \Bigg\{
  v_{r}(x,e^{-b}- e^{-t} ;b) E(e^{-b}- e^{-t},t;0,b;M)  \\
&  &
-  v_{r}(x,0 ;b) E(0,t;0,b;M)- \int_{ 0}^{ e^{-b}- e^{-t}}   v_{r}(x,r ;b) E_r(r,t;0,b;M)  \, dr \Bigg\}.
\end{eqnarray*}
The definition of $v$, (\ref{1.21}),  and (\ref{3.10}) imply
\[
\frac{1}{2} A(x,\partial_x) u(x,t)
=
 \int_{ 0}^{t} db \Bigg\{
  v_{t}(x,e^{-b}- e^{-t} ;b) (4e^{-b-t})^{-\frac{1}{2}}  
- \int_{ 0}^{ e^{-b}- e^{-t}}   v_{r}(x,r ;b) E_r(r,t;0,b;M)  \, dr \Bigg\}.
\]
Hence
\begin{eqnarray*}
&  &
\frac{1}{2} \left( \partial_t^2 u(x,t) -  e^{-2t} A(x,\partial_x) u(x,t) - M^2 u(x,t)\right) \\
&  =  &
-e^{-t}    \int_{ 0}^{t}
v(x, e^{-b}- e^{-t} ;b) (4e^{-b-t})^{-\frac{1}{2}}  \, db
+
\frac{1}{2} f(x, t)     \\
&  &
+ e^{-2t}       \int_{ 0}^{t}
v(x, e^{-b}- e^{-t} ;b) E_r( e^{-b}- e^{-t},t;0,b;M)  \, db   \\
&  &
+2e^{-t}    \int_{ 0}^{t}
v(x, e^{-b}- e^{-t} ;b) E_t( e^{-b}- e^{-t},t;0,b;M)  \, db  \\
&  &
 +e^{-2t}
 \int_{ 0}^{t} db  \int_{ 0}^{ e^{-b}- e^{-t}}   v_{r}(x,r ;b) E_r(r,t;0,b;M)  \, dr      \\
&  &
+  \int_{ 0}^{t} db
  \int_{ 0}^{ e^{-b}- e^{-t}}  v(x,r ;b)\Big\{  E_{tt} (r,t;0,b;M) -M^2 E(r,t;0,b;M) \Big\}  \, dr\,.
\end{eqnarray*}
On the other hand, the integration by parts leads to
\begin{eqnarray*}
&  &
 \int_{ 0}^{t} db  \int_{ 0}^{ e^{-b}- e^{-t}}   v_{r}(x,r ;b) E_r(r,t;0,b;M)  \, dr     \\
 & = &
  \int_{ 0}^{t} db  \Bigg\{ v(x,e^{-b}- e^{-t} ;b) E_r(e^{-b}- e^{-t},t;0,b;M) - v(x,0 ;b) E_r(0,t;0,b;M) \\
&  &
-\int_{ 0}^{ e^{-b}- e^{-t}}   v(x,r ;b) E_{rr}(r,t;0,b;M)  \, dr \Bigg\}\,.
\end{eqnarray*}
Then we use Proposition~\ref{P2.3}.
Consequently,
\begin{eqnarray*}
&  &
 \int_{ 0}^{t} db  \int_{ 0}^{ e^{-b}- e^{-t}}   v_{r}(x,r ;b) E_r(r,t;0,b;M)  \, dr     \\
 & = &
  \int_{ 0}^{t} db  \Bigg\{ v(x,e^{-b}- e^{-t} ;b) E_r(e^{-b}- e^{-t},t;0,b;M)
-\int_{ 0}^{ e^{-b}- e^{-t}}   v(x,r ;b) E_{rr}(r,t;0,b;M)  \, dr \Bigg\}\,.
\end{eqnarray*}
Thus,
\begin{eqnarray*}
&  &
\frac{1}{2} \left( \partial_t^2 u(x,t) -  e^{-2t} A(x,\partial_x) u(x,t) - M^2 u(x,t)\right) \\
&  =  &
f(x, t)  \frac{1}{2} -e^{-t}    \int_{ 0}^{t}
v(x, e^{-b}- e^{-t} ;b) (4e^{-b-t})^{-\frac{1}{2}}  \, db   \\
&  &
+ e^{-2t}       \int_{ 0}^{t}
v(x, e^{-b}- e^{-t} ;b) E_r( e^{-b}- e^{-t},t;0,b;M)  \, db   \\
&  &
+2e^{-t}    \int_{ 0}^{t}
v(x, e^{-b}- e^{-t} ;b) E_t( e^{-b}- e^{-t},t;0,b;M)  \, db  \\
&  &
 +e^{-2t}  \int_{ 0}^{t} db    v(x,e^{-b}- e^{-t} ;b) E_r(e^{-b}- e^{-t},t;0,b;M)    \\
&  &
+  \int_{ 0}^{t} db
  \int_{ 0}^{ e^{-b}- e^{-t}}  v(x,r ;b)\Big\{  E_{tt} (r,t;0,b;M) - e^{-2t} E_{rr}(r,t;0,b;M)-M^2 E(r,t;0,b;M) \Big\}  \, dr\,.
\end{eqnarray*}
Now we apply Theorem~\ref{T2.12} to the last integral and obtain
\begin{eqnarray*}
&  &
\frac{1}{2} \left( \partial_t^2 u(x,t) -  e^{-2t} A(x,\partial_x) u(x,t) - M^2 u(x,t)\right) \\
&  =  &
 \frac{1}{2}f(x, t)    +e^{-2t}   \int_{ 0}^{t}  v(x,e^{-b}- e^{-t} ;b)   \\
&  &
\times \Bigg\{ 2 E_r(e^{-b}- e^{-t},t;0,b;M) + 2e^t E_t( e^{-b}- e^{-t},t;0,b;M)- e^t (4e^{-b-t})^{-\frac{1}{2}}
 \Bigg\}    db\,.
\end{eqnarray*}
Proposition~\ref{P2.9} completes the proof of this case. It is easy to check the initial conditions. Theorem~\ref{T1.2} is proven.
\hfill $\square$

\begin{small}

\end{small}
\end{document}